\documentclass[10pt,preprint]{article}

\usepackage{amsmath, latexsym, amsfonts, amssymb, amsthm}
\usepackage{graphicx, color,hyperref,dsfont,epsfig,caption,wrapfig,subfig}
\usepackage{pstricks}
\usepackage{movie15}
\usepackage{epigraph}
\hypersetup{
    linktoc=page,
    linkcolor=red,          
    citecolor=blue,        
    filecolor=blue,      
    urlcolor=cyan,
    colorlinks=true           
}

\numberwithin{equation}{section}

\newtheorem{theorem}{Theorem}[section]
\newtheorem{lemma}[theorem]{Lemma}
\newtheorem{proposition}[theorem]{Proposition}
\newtheorem{corollary}[theorem]{Corollary}

\newtheorem{remark}[theorem]{Remark}
\newtheorem{definition}[theorem]{Definition}

\theoremstyle{definition}

\renewcommand{\tilde}{\widetilde}          
\DeclareMathSymbol{\leqslant}{\mathalpha}{AMSa}{"36} 
\DeclareMathSymbol{\geqslant}{\mathalpha}{AMSa}{"3E} 
\DeclareMathSymbol{\eset}{\mathalpha}{AMSb}{"3F}     
\renewcommand{\leq}{\;\leqslant\;}                   
\renewcommand{\geq}{\;\geqslant\;}                   

\newcommand{\C}{\mathbb{C}}

\newcommand{\R}{\mathbb{R}}
\newcommand{\Z}{\mathbb{Z}}

\newcommand{\N}{\mathbb{N}}
\newcommand{\Q}{\mathbb{Q}}

\newcommand{\E}{\mathds{E}}
\renewcommand{\P}{\mathds{P}}

\usepackage{xparse}

\DeclareDocumentCommand \Pmp { m m o} {
\IfNoValueTF{#3}
{P_{#1}^{#2}}
{P_{#1}^{#2}\left(#3\right)}
}

\DeclareDocumentCommand \Emp { m m o} {
\IfNoValueTF{#3}
{E_{#1}^{#2}}
{E_{#1}^{#2}\left[#3\right]}
}

\DeclareDocumentCommand \Pbr { m m m m o } {
\IfNoValueTF{#5}
{P_{#1}^{#2\stackrel{#4}{\rightarrow} #3}}
{P_{#1}^{#2\stackrel{#4}{\rightarrow} #3}\left(#5\right)}
}

\DeclareDocumentCommand \Ebr { m m m m o } {
\IfNoValueTF{#5}
{E_{#1}^{#2\stackrel{#4}{\rightarrow} #3}}
{E_{#1}^{#2\stackrel{#4}{\rightarrow} #3}\left[#5\right]}
}

\def\S{\mathbb{S}}

\def\bi{\begin{itemize}}
\def\ei{\end{itemize}}

\def\bnum{\begin{enumerate}}
\def\enum{\end{enumerate}}
\def\<#1{\langle #1 \rangle}

\def\caT{\mathcal{T}}

\title{Lecture notes on Liouville theory and the DOZZ formula}

 \author{ 
 Vincent Vargas \footnote{ENS Ulm, DMA, 45 rue d'Ulm,  75005 Paris, France.} }

\date{\vspace{-5ex}}

\begin{document}

\maketitle

\begin{abstract}
The purpose of these notes, based on a series of 4 lectures given by the author at IHES, is to explain the recent proof of the DOZZ formula for the three point correlation functions of Liouville conformal field theory (LCFT). We first review the probabilistic construction of the $N$ point correlation functions of LCFT on the Riemann sphere for $N \geq 3$ (based on a series of works of the author with David, Kupiainen and Rhodes). We then explain the construction of the two point correlation functions of LCFT, also called reflection coefficient. These probabilistic constructions will be justified from the point of view of Polyakov's path integral formulation of LCFT. Finally, we explain the proof of the DOZZ formula for the three point correlation functions of LCFT (based on 2 papers by the author with Kupiainen and Rhodes).
\end{abstract}
\begin{center}

\epigraph{``There are methods and {\bf formulae} in science, which
serve as master-keys to many apparently different problems."}{--- \textup{ {\bf Polyakov}}, \emph{Quantum geometry of bosonic strings}, 1981 }

\end{center}

\tableofcontents
\noindent

\section{Introduction and history of the DOZZ formula}

 In the 1981 seminal paper ``Quantum geometry of bosonic strings", Polyakov introduced a path integral based theory of summation over Riemannian metrics. In his 1981 paper, Polyakov discovered that an essential building block of the theory is Liouville conformal field theory (LCFT hereafter) which describes the conformal factor of the metric chosen ``at random". This initiated a fascinating quest in physics to understand LCFT and 2d conformal symmetry in general\footnote{Very recently, there has also been a spectacular revival in physics of CFT in 3 dimensions: see for instance the seminal works \cite{Ry1,Ry2} where CFT techniques are used to study the critical 3d Ising model.}; ultimately this led to the discovery 15 years later of the celebrated DOZZ formula \cite{DoOt,ZaZa}  which appears today in a wide variety of contexts in the physics literature and in particular recently in relation with 4d supersymmetric gauge theories (via the AGT conjecture \cite{AGT}). The purpose of these notes, based on tape recorded lectures at IHES \cite{Var}, is to present the first mathematical proof of the DOZZ formula using probability theory.
 
In Polyakov's formulation, LCFT is a 2d version of a Feynman path integral with an  {\bf exponential interaction} term. Mathematically, LCFT is an infinite measure on some functional space; we will denote the underlying random field $\phi$, the so-called Liouville field. In what follows, we give an account of LCFT at the level of physics rigour. First, we introduce the following notation for the so-called ``vertex operators"associated to $\phi$: for all $z \in \C$ and $\alpha \in \C$, we set
\begin{equation}
V_\alpha(z)= e^{\alpha \phi(z)}.
\end{equation} 

The purpose of LCFT on the Riemann sphere $\S^2= \C \cup \lbrace \infty \rbrace$ is to study the fields $V_\alpha(z)$ under a natural measure (defined by the so-called Liouville action) and this amounts to computing the following correlations for all distinct $(z_k)_{1 \leq k \leq N} \in \C^N$ and $(\alpha_k)_{1 \leq k \leq N} \in \C^N$: 
\begin{equation}\label{defcorrelations}
\langle \prod_{k=1}^N V_{\alpha_k} (z_k) \rangle_{\gamma,\mu}  := \int_{\varphi: \S^2 \to \R} \left ( \prod_{k=1}^N e^{ \alpha_k \phi(z_k) }\right ) e^{-  S_L(\varphi) } D \varphi, 
\end{equation}
where: 
\begin{itemize}

\item
$D\varphi$ is the ``Lebesgue measure" on the space of real functions $\varphi :\S^2 \to \R$.  
\item
 $S_L$ is the Liouville action: 
\begin{equation}\label{defLiouvilleaction}
S_L(\varphi):= \frac{1}{4\pi}\int_{\mathbb{S}^2}\big(  |\nabla_g \varphi(x) |^2 +R_g(x) Q \varphi(x)  +4\pi \mu e^{\gamma \varphi(x)  } \big) g(x)\, d^2 x 
\end{equation}
where $d^2x$ denotes the standard Lebesgue measure and:
\begin{itemize}
 \item
 $g(x)|dx|^2$ is some metric on $\S^2$  (with associated gradient $\nabla_g$)
 \item
 $R_g(x)=-\frac{1}{g(x)}\Delta \ln g(x)$ is the curvature  of $g$ ($\Delta$ denotes the standard Laplacian on $\C$)
  \item
 $\gamma$ is a positive parameter belonging to  $(0,2)$ and $Q=\frac{\gamma}{2}+\frac{2}{\gamma}$
\item
$\mu>0$ is the so-called cosmological constant
\end{itemize}

\item
$\phi(x)= \varphi(x)+\frac{Q}{2} \ln g(x)$ is the {\bf Liouville field}.

\end{itemize}

Since the Lebesgue measure does not exist on an infinite dimensional space, one can see that the correlations \eqref{defcorrelations} are ill-defined mathematically. One of the purpose of these notes is to explain how to actually give a probabilistic meaning to \eqref{defcorrelations}: this requires some care. Moreover, we will see that in the probabilistic definition of the measure $e^{-  S_L(\varphi) } D \varphi$, the regularity of $\phi$ is that of a distribution (in the sense of Schwartz) so actually making sense of the correlations of $V_\alpha(z)$ or the interaction term $\mu \int_{\mathbb{S}^2} e^{\gamma \varphi(x)  }  g(x)\, d^2 x $ requires a renormalization procedure.  

LCFT is a conformal field theory (CFT) which implies that the correlations  \eqref{defcorrelations} behave  like conformal tensors: see \eqref{KPZformula} below for a precise statement.   These conformal properties will be instrumental in stating a precise conjecture relating the scaling limit of planar maps to LCFT: this is the content of the celebrated KPZ conjecture (Lecture 3).

Now, we turn to the computation of the correlations in the physics literature. Following Polyakov's 1981 work, Belavin-Polyakov-Zamolodchikov (BPZ) laid in their 1984 seminal work \cite{BPZ} the foundations of modern CFT. Their aim was to give explicit expressions for \eqref{defcorrelations} and they proposed a general scheme to compute correlations for CFTs. Roughly, they proposed formulae based on the knowledge of the 3 point correlation functions and a recursive procedure (the so-called conformal bootstrap) to compute the higher order correlations. Unfortunately, they could not apply their powerful methodology to LCFT because they were unable to propose an exact expression for the three point correlation functions of the theory.

In the middle of the 90's, an intriguing formula for the three point correlation functions was proposed independently by Dorn-Otto \cite{DoOt} and Zamolodchikov-Zamolodchikov (DOZZ) \cite{ZaZaarxiv,ZaZa}, the so-called {\bf DOZZ formula}. We will see at the end of lecture 1 their derivation of the DOZZ formula: it is based on an ingenious analytic continuation ``guess". If we set
\begin{equation}
C_\gamma(\alpha_1,\alpha_2,\alpha_3) := \underset{|z| \to \infty}{\lim} |z|^{2\alpha_3(Q-\frac{\alpha_3}{2})} \langle V_{\alpha_1} (0)V_{\alpha_2} (1)V_{\alpha_3} (z) \rangle_{\gamma,\mu}\footnote{We will see that conformal symmetry explains the renormalization term $ |z|^{2\alpha_3(Q-\frac{\alpha_3}{2})}$. $C_\gamma (\alpha_1,\alpha_2,\alpha_3)$ is often denoted quite naturally $\langle V_{\alpha_1} (0)V_{\alpha_2} (1)V_{\alpha_3} (\infty) \rangle_{\gamma,\mu}$ in the physics literature.}
\end{equation}
and $l(z)=\frac{\Gamma (z)}{\Gamma (1-z)}$ then DOZZ argued that the following exact formula should hold for $\alpha_1,\alpha_2,\alpha_3 \in \C$ with $\bar{\alpha}=\alpha_1+\alpha_2+\alpha_3$ 
\begin{equation}\label{theDOZZformula}
C_\gamma(\alpha_1,\alpha_2,\alpha_3)  = (\pi \: \mu \:  l(\frac{\gamma^2}{4})  \: (\frac{\gamma}{2})^{2 -\gamma^2/2} )^{\frac{2 Q -\bar{\alpha}}{\gamma}}   \frac{\Upsilon_{\frac{\gamma}{2}}'(0)\Upsilon_{\frac{\gamma}{2}}(\alpha_1) \Upsilon_{\frac{\gamma}{2}}(\alpha_2) \Upsilon_{\frac{\gamma}{2}}(\alpha_3)}{\Upsilon_{\frac{\gamma}{2}}(\frac{\bar{\alpha}}{2}-Q) 
\Upsilon_{\frac{\gamma}{2}}(\frac{\bar{\alpha}}{2}-\alpha_1) \Upsilon_{\frac{\gamma}{2}}(\frac{\bar{\alpha}}{2}-\alpha_2) \Upsilon_{\frac{\gamma}{2}}(\frac{\bar{\alpha}}{2} -\alpha_3)   }
\end{equation}      
where  $\Upsilon_{\frac{\gamma}{2}}(z)$ is Zamolodchikov's special holomorphic function defined on $\C$ with the following 
explicit expression for $0<\Re (z)< Q$\footnote{The function 
has a simple construction in terms of standard double gamma functions or the double zeta function: see the reviews \cite{nakayama,Rib} or the paper \cite{Tesc1} for instance.}
\begin{equation}\label{def:upsilon}
\ln \Upsilon_{\frac{\gamma}{2}} (z)  = \int_{0}^\infty  \left ( \Big (\frac{Q}{2}-z \Big )^2  e^{-t}-  \frac{( \sinh( (\frac{Q}{2}-z )\frac{t}{2}  )   )^2}{\sinh (\frac{t \gamma}{4}) \sinh( \frac{t}{\gamma} )}    \right ) \frac{dt}{t}.
\end{equation}
We will denote $C_\gamma^{{\rm DOZZ}}(\alpha_1,\alpha_2,\alpha_3)$ the DOZZ formula, i.e. the right hand side of \eqref{theDOZZformula}. In fact, the function $\Upsilon_{\frac{\gamma}{2}}$  is defined on $\C$ by analytic continuation of the expression \eqref{def:upsilon} as expression \eqref{def:upsilon}  satisfies the following remarkable functional relations: 
\begin{equation}\label{shiftUpsilon}
\Upsilon_{\frac{\gamma}{2}} (z+\frac{\gamma}{2}) = l( \frac{\gamma}{2}z) (\frac{\gamma}{2})^{1-\gamma z}\Upsilon_{\frac{\gamma}{2}} (z), \quad
\Upsilon_{\frac{\gamma}{2}} (z+\frac{2}{\gamma}) = l(\frac{2}{\gamma}z) (\frac{\gamma}{2})^{\frac{4}{\gamma} z-1} \Upsilon_{\frac{\gamma}{2}} (z).
\end{equation}
The function $\Upsilon_{\frac{\gamma}{2}}$ has no poles in $\C$ and the zeros of $\Upsilon_{\frac{\gamma}{2}}$ are simple (if $\gamma^2 \not \in \Q$) and given by the discrete set $(-\frac{\gamma}{2} \N-\frac{2}{\gamma} \N) \cup (Q+\frac{\gamma}{2} \N+\frac{2}{\gamma} \N )$.

The purpose of these notes is precisely to present a mathematical proof of the DOZZ formula \eqref{theDOZZformula}. Let us stress that identity  \eqref{theDOZZformula} is far from obvious and even controversial in the physics literature. This is due to the fact that the DOZZ formula  $C_\gamma^{{\rm DOZZ}}(\alpha_1,\alpha_2,\alpha_3)$  is invariant under the substitution of parameters
 \begin{equation}\label{dualityDOZZ}
 \frac{\gamma}{2} \leftrightarrow \frac{2}{\gamma}, \quad  \mu  \leftrightarrow \tilde{\mu}= \frac{(\mu \pi \ell(\frac{\gamma^2}{4})  )^{\frac{4}{\gamma^2} }}{ \pi \ell(\frac{4}{\gamma^2})}.
\end{equation}
This  duality symmetry is at the core of the DOZZ controversy. Indeed this symmetry is not manifest in the Liouville action functional \eqref{defLiouvilleaction} and therefore duality is axiomatically assumed in the physics literature. As a matter of fact, it has been argued in physics that it could be necessary to add the {\bf dual} potential $4 \pi \tilde{\mu}e^{\frac{4}{\gamma} \varphi(x)}$ with $\tilde{\mu}$ defined in \eqref{dualityDOZZ} in the  Liouville action functional \eqref{defLiouvilleaction} for the DOZZ formula to hold (see \cite{OPS, Tesc1, ZaZaarxiv, ZaZa}). The proof of the DOZZ formula presented in these notes is based on the Liouville action functional \eqref{defLiouvilleaction} and therefore demonstrates that adding a dual potential is not necessary. 

These notes are divided into 4 parts which are meant to be read in this order:
\begin{itemize}
\item
Lecture 1: in this lecture, we give a rigorous probabilistic definition of the correlations \eqref{defcorrelations}. The probabilistic construction will be based on the Gaussian Free Field (GFF hereafter) and the exponential of the GFF, i.e. Kahane's Gaussian multiplicative chaos (GMC hereafter). We will in particular define mathematically \eqref{defcorrelations} for real $\alpha_i$ satisfying certain bounds: this will enable us to give a mathematical content to the DOZZ formula \eqref{theDOZZformula}. In lecture 1, we will not explain in detail why our probabilistic definitions are a faithful definition of \eqref{defcorrelations} based on the formal path integral: this will be the topic of lecture 3.

\item
Lecture 2: in this lecture, we will explain how to construct in a probabilistic way the two point correlation functions of LCFT, also called reflection coefficient in LCFT. The construction is based on a limiting procedure on certain three point correlation functions. These two point correlation functions admit an explicit expression (as a corollary of the DOZZ formula). We will see that the key ingredient is to understand the tail expansion of GMC measures. We will also recall the definition of a random measure called the quantum sphere (with two marked points $0$ and $\infty$) introduced by Duplantier-Miller-Sheffield \cite{DMS} in an approach to LCFT based on a GFF/SLE coupling. We will argue that their theory corresponds to these two point correlation functions; more specifically, the two point correlation functions can be seen as the partition function of their theory.   

\item
Lecture 3: we will justify in this lecture the probabilistic definitions of lecture 1. More specifically, we will show how an appropriate interpretation of the path integral leads to the definitions of lecture 1. We will also prove certain conformal invariance properties of LCFT; finally, we explain how these conformal invariance properties lead naturally to the KPZ conjecture on random planar maps.

\item
Lecture 4: we oultine in this lecture the proof of the DOZZ formula. The proof is based on the BPZ differential equations and also on rigorous operator product expansions (OPE hereafter); in the context of LCFT, an OPE is nothing but a Taylor expansion when two points of the correlation functions collide. The GMC tail expansions and the reflection coefficient enter the OPE in an essential way here.

\end{itemize}

\vspace{0.3 cm}

\subsubsection*{Some historical references and reviews from physics}
In the physics literature, there are numerous well-written reviews on LCFT. Most of these reviews contain an interesting discussion on the path integral construction of LCFT but essentially develop another approach to LCFT, the conformal bootstrap approach. In these lectures, we will not discuss the conformal bootstrap approach which is based on two ingredients: in this approach applied to LCFT, the three point correlation functions are set axiomatically to be defined (up to conformal invariance) by the DOZZ formula and the higher order correlation functions are constructed using a recursive procedure\footnote{Actually, the proof of the DOZZ formula presented in these lectures is part of a program which aims at reconciling the path integral approach and the conformal bootstrap approach to LCFT.}. Let us mention and comment (in chronological order) the following reviews which we find enlightening:
\begin{itemize}
\item
The 1990 review by Seiberg \cite{seiberg} ``Notes on Quantum Liouville Theory and Quantum Gravity". This review contains many pioneering ideas on the path integral formulation of LCFT (like the existence of some form of two point correlation for LCFT: see lecture 2 for a rigorous point of view).   
\item
The 2001 paper by Teschner \cite{Tesc1} ``Liouville theory revisited" which develops the so-called vertex operator approach to LCFT. The paper  \cite{Tesc1} is not just a review as it proposes strong arguments to justify the use of the DOZZ formula as a building block in the conformal bootstrap approach.
\item
The 2004 review by Nakayama \cite{nakayama} ``Liouville Field Theory --A decade after the revolution  " which contains a very complete and exhaustive account on LCFT and quantum gravity. 
\item
The 2014 review by Ribault \cite{Rib} ``Conformal field theory on the plane" which exposes in a self consistent way the conformal bootstrap approach to LCFT and more generally of many CFTs. The CFT point of view exposed in Ribault is axiomatic and does not rely on any underlying path integral. The methodology is powerful and in particular the review contains an account on the very recent construction by Ribault-Santachiara of LCFT with imaginary $\gamma$ \cite{Rib1} (which is based on an intriguing ``imaginary" DOZZ formula discovered by Schomerus \cite{Scho}, Kostov-Petkova \cite{KoPe} and Zamolodchikov \cite{Za} ). It is an open problem to give a probabilistic and more generally mathematical content to ``imaginary" LCFT.   
\end{itemize}

Finally, on the mathematical side, there has been an attempt to make sense of the path integral formulation of LCFT  by Takhtajan-Teo \cite{TT} using  the formalism of classical conformal geometry. In this formalism, LCFT was defined as a formal power series in $\gamma$. On the probabilistic side, numerous works have attempted to give a rigorous definition to the path integral using GFF theory but the work of David-Kupiainen-Rhodes-Vargas \cite{DKRV} is the first to give a complete and rigorous probabilistic defintion to \eqref{defcorrelations} (in physics, there is also a GFF approach to the path integral formulation of LCFT which was first developed in Goulian-Li \cite{GouLi}).

\subsubsection*{Conventions and notations} In what follows,  $z$,  $x,y$ and   $z_1,\dots, z_N$ all denote complex variables. We use the standard notation for complex derivatives $\partial_x= \frac{1}{2}(\partial_{x_1}- i \partial_{x_2})$ and $\partial_{\bar{x}}= \frac{1}{2} (\partial_{x_1}+ i \partial_{x_2})$ for $x=x_1+ix_2$. The  Lebesgue measure on $\C$ (seen as $\R^2$) is denoted by $d^2x$, the standard gradient  is denoted $\nabla$ and the standard Laplacian is denoted $\Delta$. We will also denote $|\cdot|$ the norm in $\C$ of the standard Euclidean (flat) metric.

 Let $\mathcal{S}(\C)$ be the standard space of smooth functions $f$ with fast decay at infinity, i.e. the functions $f$ that are $C^\infty$ and such that any derivative decays at infinity faster than any inverse polynomial.  We denote by $\mathcal{S}'(\C)$ the associated space of tempered distributions. We also set $\mathcal{S}_0(\C)$ the subset of functions with vanishing average with respect to the Lebesgue measure $d^2x$, i.e. such that $\int_\C f(x) d^2x=0$; the dual of $\mathcal{S}_0(\C)$ is then the space $\mathcal{S}'(\C) / \mathbb{R}$ of distributions defined up to a global constant. If $f$ belongs to $\mathcal{S}(\C)$ then we denote $\hat{f}(\xi)= \int_{\C}   e^{2i \pi \xi.x} f(x) d^2x$ its Fourier transform.

Let us recall here the main notations in these notes (exact definitions will appear in the lectures). In these notes, we will define and work with:
\begin{itemize}
\item
The special function $l$ defined by $l(x)=\frac{\Gamma(x)}{\Gamma (1-x)}$ where $\Gamma$ is the standard Gamma function.
\item
The cosmological constant $\mu>0$ and the {\bf dual} cosmological constant $\tilde{\mu}=\frac{(\mu \pi l(\frac{\gamma^2}{2})  )^{\frac{4}{\gamma^2} }}{ \pi l(\frac{4}{\gamma^2})}$.
\item
A parameter $\gamma \in (0,2)$ and a parameter $Q$ satisfying $Q=\frac{\gamma}{2}+\frac{2}{\gamma}$.
\item
A conformal metric $g$ given by $g(x)=\frac{1}{|x|_+^4}$ where $|x|_+=|x| \vee 1 $ ($\vee$ denotes the maximum). In this setting, the volume with respect to $g$ is simply $g(x)d^2x$. The gradient with respect to $g$ is $\Delta_g= \frac{1}{g} \nabla$ and the Laplacian is $\Delta_g= \frac{1}{g}\Delta$.  The curvature of $g$ defined by $-\frac{1}{g(x)}\Delta \ln g(x)$ is a Radon measure and more precisely 4 times the uniform measure on the unit circle. 
\item
The Gaussian Free Field (GFF) $X$ with vanishing mean on the unit circle or equivalently with vanishing mean along the curvature of $g$. The covariance of $X$ is $\E[X(x)X(y)]= \ln \frac{1}{|x-y|}+\ln |x|_++\ln |y|_+$. For all $\epsilon>0$, $X_\epsilon(x)$ is a regularization of $X$ at scale $\epsilon$; in most of these notes, $X_\epsilon(x)$ will denote the circle average of $X$ on a circle of center $x$ and radius $\epsilon$ (with respect to the standard Euclidean distance).
\item
The exponential of $X$, i.e. a Gaussian multiplicative chaos (GMC) measure defined formally by the formula $M_\gamma(d^2x)= e^{\gamma X(x)-\frac{\gamma^2}{2} \E[X(x)^2]  } g(x) d^2x$. We will use the notation $M_\gamma(d^2x)$ as well as $e^{\gamma X(x)-\frac{\gamma^2}{2} \E[X(x)^2]  } g(x) d^2x$. 
\item
For all $\alpha \in (\frac{\gamma}{2},Q)$, $(\mathcal{B}^\alpha_s)_{s \in \R}$ denotes a two sided Brownian motion starting from $0$ at $s=0$ with negative drift $\alpha-Q$ on each side of $0$ and conditioned to be negative.
\item
The Gaussian lateral noise field $Y(t,\theta)$ defined on the cylinder $\R \times [0,2\pi]$ with covariance $\E[  Y(s,\theta) Y(t,\theta') ]  = \ln \frac{e^{-s}\vee e^{-t}}{|e^{-s}e^{i \theta} - e^{-t} e^{i \theta'} |}$. We will also consider the GMC measure $N_\gamma(ds d\theta):=e^{\gamma Y(s,\theta)-\frac{ \gamma^2 }{2}E[Y(s,\theta)^2]}  ds d\theta$ with respect to the field $Y$.   We will especially be interested in the (stationary) total mass process
$Z_s = \int_0^{2 \pi} e^{\gamma Y(s,\theta)-\frac{\gamma^2 }{2}E[Y(s,\theta)^2]} d\theta$.

\end{itemize}

\subsubsection*{Acknowledgements} The author would first like to thank Hugo Duminil-Copin for inviting him to give lectures at IHES. He would also like to thank those who read a preliminary draft of these notes: Yichao Huang, Guillaume Remy, R\'emi Rhodes, Sylvain Ribault, Tunan Zhu. Their valuable comments have certainly improved the readability of these notes.

\section{Lecture 1: probabilistic statement of the DOZZ formula}

In the first lecture, we introduce the $N$ point correlation functions of LCFT for $N \geq 3$. The probabilistic construction is based on the GFF and GMC (formally the exponential of the GFF) so we review these two essential building blocks of the theory (one can also have a look at the reviews by Rhodes-Vargas \cite{RV} or Kupiainen \cite{cargese}). With this material, we can then  introduce a rigorous mathematical definition of \eqref{defcorrelations} and give a precise mathematical content to the DOZZ formula. The fact that this rigorous mathematical definition is a faithful definition of  \eqref{defcorrelations} will be explained in lecture 3. Finally, we explain how DOZZ derived the DOZZ formula in the middle of the 90's.
 
\subsection{Gaussian Free Field and Gaussian multiplicative chaos theory}

\subsubsection{The Gaussian Free Field}
We introduce the GFF we will be working with in these notes (we refer to section 4 in Dub\'edat's paper \cite{Dub0} or the review by Sheffield \cite{She07} for an introduction to the GFF).  If $T$ is some random field living in $\mathcal{S}'(\C)$ or $\mathcal{S}'(\C) / \mathbb{R}$,  we adopt the notation of generalized functions by denoting $\int_\C T(x) f(x)d^2x$ the field $T$ evaluated at $f \in \mathcal{S}(\C)$ or $f \in \mathcal{S}_0(\C)$.     

We start by introducing the full plane GFF $\bar{X}$ as as a Gaussian random variable living in $\mathcal{S}'(\C) / \mathbb{R}$:

\begin{definition}
 The full plane GFF $\bar{X}$ is defined as a Gaussian random variable living in $\mathcal{S}'(\C) / \mathbb{R}$ with the following covariance for all $f,h \in \mathcal{S}_0(\C)$
 \begin{equation}\label{defcorrelsbarX}
 \E\left [ \left (\int_\C f(x) \bar{X}(x) d^2x \right )  \left ( \int_\C h(x) \bar{X}(x) d^2x  \right ) \right   ]=   \int_{\C^2}  f(x) h(y)  \ln \frac{1}{|x-y|} d^2x d^2y
 \end{equation} 
 \end{definition}
The above definition defines a Gaussian field since one can prove by using Fourier analysis that for all $f,h \in \mathcal{S}_0(\C)$ 
 \begin{equation}
  \int_{\C^2}  f(x) h(y)  \ln \frac{1}{|x-y|} d^2x d^2y= c \int_{\C}  \frac{\hat{f}(\xi) \overline{\hat{h}(\xi)}}{|\xi|^2}d^2\xi
 \end{equation} 
where $c>0$ is some constant. One can notice that the right hand side of \eqref{defcorrelsbarX} is well defined for very irregular functions $f,g$ and even for generalized functions. Therefore, one can consider the variable $\int_\C f(x) \bar{X}(x) d^2x$  for general $f$ and we will do so in the sequel with no further notice.

Once the full plane GFF is introduced, one can construct different GFFs living in $ \mathcal{S}'(\C)$ (hence defined against any test function in $\mathcal{S}(\C)$) by prescribing the value of $\bar{X}$ against a probability measure $\rho(d^2x)$. More specifically, if $\rho(d^2x)$ is a probability measure such that $\int_{\C^2}   \ln \frac{1}{|x-y|} \rho(d^2x) \rho(d^2y)< \infty$, one can introduce the GFF $X_\rho$ with average $0$ on $\rho(d^2x)$ by the following formula
\begin{equation}
X_\rho(x)= \bar{X}(x)- \int_{\C} \bar{X}(x) \rho (d^2x)
\end{equation}  
Let us review two fundamental examples. If $\rho(d^2x)= \frac{4}{(1+|x|^2)^2} d^2x$ then $X_\rho$ is the GFF with vanishing mean on the Riemann sphere. If $\rho$ is the uniform probability measure on the unit circle then $X_\rho$ is the GFF with vanishing mean on the unit circle. We will denote the latter $X$ as we will only work with this GFF in the rest of these notes. A rather straightforward computation yields for all $f,h \in \mathcal{S}(\C)$
\begin{equation}\label{CovXprecise}
\E\left [ \left (\int_\C f(x) X(x) d^2x \right )  \left ( \int_\C h(x) X(x) d^2x  \right ) \right   ]= \int_{\C^2}  f(x) h(y)  K(x,y) d^2x d^2y 
\end{equation}
where $K(x,y)=\ln \frac{1}{|x-y|}+ \ln |x|_++\ln |y|_+ $ with $|x|_+= \max(  |x|,1 )$. In the sequel, we will usually write (with an obvious abuse of notation) that 
\begin{equation}\label{CovX}
\E[X(x)X(y)]=\ln \frac{1}{|x-y|}+ \ln |x|_++\ln |y|_+. 
\end{equation}

\subsubsection{Gaussian multiplicative chaos theory}

Now we introduce the GMC measure associated to $X$, i.e. the random measure $M_\gamma (d^2x)$ with (formal) density $e^{\gamma X(x)-\frac{\gamma^2}{2}\E[X(x)^2]}$ with respect to the measure $g(x) d^2x=\frac{d^2x}{|x|_+^4}$. Since $X(x)$ is not defined pointwise, the rigorous construction of $M_\gamma(d^2x)$ requires a renormalization procedure. More precisely, it is natural to construct $M_\gamma(d^2x)$ as the limit when $\epsilon$ goes to $0$ of $M_{\epsilon,\gamma}(d^2x)$ defined by  
\begin{equation}\label{associX}
M_{\epsilon,\gamma}(d^2x)= e^{\gamma X_\epsilon(x)-\frac{\gamma^2 }{2}\E[X_\epsilon(x)^2]} g(x)d^2x
\end{equation}
where $(X_\epsilon)_{\epsilon>0}$ is a sequence of smooth continuous fields which approximate $X$. We will not try to give a general definition of smooth continuous fields which approximate $X$: let us just mention that if $\theta$ is some smooth function of average $1$ then  $X_\epsilon(x)=\int_{\C}  \frac{1}{\epsilon^2} \theta (\frac{x-y}{\epsilon})X(y)d^2y$ is a sequence of smooth continuous fields which approximate $X$. As a matter of fact, one can even choose $\theta(y) d^2y$ to be the uniform measure on the unit circle (with respect to the standard Euclidean distance) in which case $X_\epsilon(x)$ is the average of the field $X$ on the circle of center $x$ and radius $\epsilon$. Let us now state a precise proposition: 

\begin{proposition}[Kahane, 1985]
Let $\gamma \in (0,2)$ and $(X_\epsilon)_{\epsilon>0}$ be a sequence of smooth continuous fields which approximate $X$ and  with associated $M_{\epsilon,\gamma}(d^2x)$ given by \eqref{associX}. Then one gets the following convergence in probability in the space of Radon measures 
\begin{equation}\label{defMgamma}
M_\gamma(d^2x):= \underset{\epsilon \to 0}{\lim} \: M_{\epsilon,\gamma}(d^2x)
\end{equation}
The random measure $M_\gamma(d^2x)$, which does not depend on the approximation $(X_\epsilon)_{\epsilon>0}$, is the GMC measure associated to $X$. The measure $M_\gamma(d^2x)$ is different from $0$.
\end{proposition}

\begin{remark}
In fact, Kahane proved the convergence \eqref{defMgamma} in a more restrictive setting than the one we mention here: for the latest results, we refer to Berestycki \cite{Ber} for an elementary approach and references. Let us also mention that one can also construct the measure $M_\gamma(d^2x)$  in the setting of a discrete GFF: in this setting, the field $X_\epsilon$ is a discrete GFF defined on a lattice with mesh size $\epsilon$ and  the convergence of \eqref{defMgamma} is replaced by a convergence in law (see the review by Rhodes-Vargas \cite{review} for precise statements). 
\end{remark}

\begin{remark}
One can show that if $\gamma \geq 2$ then $M_{\epsilon,\gamma}(d^2x)$ defined by \eqref{associX} converges to $0$ as $\epsilon$ goes to $0$; this is why we suppose that $\gamma \in (0,2)$ in these notes.
\end{remark}
In the sequel, $X_\epsilon$ will always denote the circle average regularization of $X$. Now, let us recall a few facts on the moments of $M_\gamma(d^2x)$. The first one is classical in the field of GMC; for any (non empty) open subset $\mathcal{O}\subset \C$ one has the following equivalence (see the reviews by Rhodes-Vargas \cite{review, RV} for instance):
\begin{equation}\label{Momexis}
\E[  M_\gamma (\mathcal{O})^p  ]< + \infty \: \: \Longleftrightarrow  \: \:  p< \frac{4}{\gamma^2}.
\end{equation}

In fact, we will need the following strengthening of \eqref{Momexis}: if $\alpha \in \R$ and $z\in \C$ then  
\begin{equation}\label{Momexisfort}
\E \left [ \left  (   \int_{|x-z| \leq 1}  \frac{1}{|x-z|^{\alpha \gamma}}   M_\gamma (d^2x)      \right )^p \right  ]< + \infty \: \: \Longleftrightarrow  \: \:  p< \frac{4}{\gamma^2} \wedge \frac{2}{\gamma} (Q-\alpha).
\end{equation}
where recall that $Q=\frac{\gamma}{2}+\frac{2}{\gamma}$. Notice that the singularity $\alpha$ has an effect on the existence of moments (compared to the case $\alpha=0$) if and only if $\alpha> \frac{\gamma}{2}$ (since $\frac{2}{\gamma} (Q-\frac{\gamma}{2})= \frac{4}{\gamma^2} $). For $\alpha \geq Q$, the singularity $ \frac{1}{|x-z|^{\alpha \gamma}}  $ is no longer integrable by $M_\gamma (d^2x)  $ and one has almost surely  
\begin{equation}\label{alphagreaterQ}
 \int_{|x-z| \leq 1}  \frac{1}{|x-z|^{\alpha \gamma}}   M_\gamma (d^2x)    = \infty. 
\end{equation}
We will admit the statements \eqref{Momexis}, \eqref{Momexisfort}, \eqref{alphagreaterQ} and refer to lecture 2 for a proof. In fact, we will prove much stronger statements in lecture 2 than \eqref{Momexisfort}  since we will give precise tail expansions for $\int_{|x-z| \leq 1}  \frac{1}{|x-z|^{\alpha \gamma}}   M_\gamma (d^2x)   $.

\subsection{The correlation functions of LCFT and the DOZZ formula}
Now, we introduce the probabilistic construction of the LCFT correlation functions and state a precise theorem on the three point correlation functions. For a justification from the path integral perspective, we refer to lecture 3. Let $\alpha_1, \cdots \alpha_N$ be $N$ real numbers and $z_1,\cdots, z_N$ distinct points in the complex plane $\C$. We introduce
\begin{equation}\label{defs}
s= \frac{\sum_{k=1}^N \alpha_k-2Q}{\gamma}
\end{equation}
and we suppose that the following condition holds
\begin{equation}\label{ThextendedSeibergbounds}
-s < \frac{4}{\gamma^2} \wedge \min_{1 \leq k \leq N}  \frac{2}{\gamma}(Q-\alpha_k), \quad \quad \alpha_k<Q, \; \; \forall k
\end{equation}
Notice that condition \eqref{ThextendedSeibergbounds} implies $N \geq 3$. Then we take as definition of the N point correlation function \eqref{defcorrelations} the following probabilistic formula (first introduced in David-Kupiainen-Rhodes-Vargas \cite{DKRV})

\begin{equation}\label{Firstdefcorrel}
 \langle    \prod_{k=1}^N V_{\alpha_k}(z_k)  \rangle_{\gamma,\mu} =2 \mu^{-s} \gamma^{-1}\Gamma(s)
 \prod_{i < j} \frac{1}{|z_i-z_j|^{\alpha_i \alpha_j}}\E \left [  \left (  \int_{\C}  F(x,{\bf z}) M_\gamma(d^2x)  \right )^{-s}  \right ]\footnote{The expectation term $\E[.]$ is a manifestation of the Liouville potential $e^{\gamma \varphi(x)}$ in the Liouville action \eqref{defLiouvilleaction}: see lecture 3 for more on this.}
\end{equation}
where 
\begin{equation}\label{coulomb}
F(x,{\bf z})=\prod_{k=1}^N \left ( \frac{ |x|_+}{|x-z_k|}  \right )^{\gamma \alpha_k} .
\end{equation}
It is a simple consequence of  \eqref{Momexis}, \eqref{Momexisfort} and \eqref{alphagreaterQ} that the expectation $\E[.]$ in  \eqref{Firstdefcorrel} is non trivial, i.e. belongs to $(0,\infty)$ if and only if \eqref{ThextendedSeibergbounds} holds. Notice that \eqref{Firstdefcorrel} has poles for $s=-n$ with $n \in \mathbb{N}$  which correspond to the poles of the  $\Gamma$ function\footnote{the function $\Gamma(x)$ has simple poles at $x=-n$ with residue $\frac{(-1)^n}{n!}$.}. 

\begin{remark}
Condition \eqref{ThextendedSeibergbounds} for existence of $\langle  \prod_{k=1}^N V_{\alpha_k}(z_k)  \rangle_{\gamma,\mu}$ (in the path integral formulation) extends the so-called Seiberg bounds \cite{seiberg} in the physics literature; the Seiberg bounds are the following conditions
\begin{equation}\label{TheSeibergbounds}
0<s, \quad \quad \alpha_k<Q, \; \; \forall k
\end{equation} 
\end{remark}

LCFT is a CFT and one can show the following KPZ type relation (after Knizhnik-Polyakov-Zamolodchikov \cite{cf:KPZ}) when $\psi(z)= \frac{az+b}{cz+d}$ is a M\"obius transform (see lecture 3 for a sketch of proof)

 \begin{equation}\label{KPZformula}
\langle \prod_{k=1}^N V_{\alpha_k}(\psi(z_k))    \rangle_{\gamma,\mu}=  \prod_{k=1}^N |\psi'(z_k)|^{-2 \Delta_{\alpha_k}}     \langle \prod_{k=1}^N V_{\alpha_k}(z_k)     \rangle_{\gamma,\mu}
\end{equation}  
where  $\Delta_{\alpha}=\frac{\alpha}{2}(Q-\frac{\alpha}{2})$ is called the conformal weight. This global conformal symmetry fixes the three point correlation functions up to a constant: 
\begin{equation}\label{firstdefstructure}
 \langle      \prod_{k=1}^3 V_{\alpha_k}(z_k)  \rangle_{\gamma,\mu} 
 =  |z_1-z_2|^{ 2 \Delta_{12}}  |z_2-z_3|^{ 2 \Delta_{23}} |z_1-z_3|^{ 2 \Delta_{13}}C_\gamma(\alpha_1,\alpha_2,\alpha_3)
\end{equation}
with $\Delta_{12}= \Delta_{\alpha_3}-\Delta_{\alpha_1}-\Delta_{\alpha_2}$, etc... The constants $C_\gamma(\alpha_1, \alpha_2, \alpha_3)$ are called the three point structure constants and they can be recovered by the following limit
\begin{equation}\label{Climit}
C_\gamma(\alpha_1,\alpha_2,\alpha_3)=\lim_{z_3\to\infty} |z_3|^{4 \Delta_3} \langle       V_{\alpha_1}(0)  V_{\alpha_2}(1) V_{\alpha_3}(z_3) \rangle.
\end{equation}
Combining the above considerations, we get 
\begin{equation}\label{expression3pointstruct}
C_\gamma(\alpha_1,\alpha_2,\alpha_3)=
2 \mu^{-s}  \gamma^{-1}\Gamma(s)  \E (\rho(\alpha_1,\alpha_2,\alpha_3)^{-s})
\end{equation}
where
\begin{equation}\label{defrho3}
\rho(\alpha_1,\alpha_2,\alpha_3)= 
 \int_{\C} 
  \frac{|x|_+^{\gamma(\alpha_1+\alpha_2+\alpha_3)}}{ |x|^{\gamma \alpha_1}  |x-1|^{\gamma \alpha_2}  }  
  M_\gamma(d^2x).
\end{equation}

We can now state our probabilistic formulation of the DOZZ formula:

\begin{theorem}[Kupiainen, Rhodes, Vargas]\label{theoremDOZZ}
Let $\alpha_1,\alpha_2,\alpha_3$ satisfy the bounds \eqref{ThextendedSeibergbounds} with $N=3$. The following equality  holds $$C_\gamma(\alpha_1,\alpha_2,\alpha_3)=C_\gamma^{{\rm DOZZ}}(\alpha_1,\alpha_2,\alpha_3).$$
\end{theorem}

From the purely probabilistic point of view, Theorem \ref{theoremDOZZ} can be interpreted as a far reaching integrability result on GMC on the Riemann sphere; indeed recall that $C_\gamma(\alpha_1,\alpha_2,\alpha_3)$ has an expression in terms of a fractional moment of some form of GMC: see formula \eqref{expression3pointstruct}. There are numerous integrability results on GMC in the physics literature; to the best of our knowledge, Theorem \ref{theoremDOZZ} is the first rigorous non trivial integrability result on GMC; we believe the techniques developped in \cite{KRV,KRV1} to prove Theorem \ref{theoremDOZZ} can be adapted to other contexts to prove many other integrability results for GMC. For instance, let us mention the recent proof by Remy \cite{remy} of the Fyodorov-Bouchaud formula \cite{fybu} for GMC on the unit circle; his proof is done in the context of LCFT in the unit disk (defined by Huang-Rhodes-Vargas \cite{HRV}) where he shows that proving he Fyodorov-Bouchaud formula is equivalent to computing certain one point correlation functions of LCFT in the unit disk.

\subsection{The derivation of the DOZZ  formula in physics: analytic continuation of the Dotsenko-Fateev integrals}

Here, we explain how Dorn-Otto-Zamolodchikov-Zamolodchikov derived their celebrated DOZZ formula from the path integral. In fact, we can actually translate their derivation within our probabilistic setup. Though formula \eqref{expression3pointstruct} was unknown to physicists for general $\alpha_1,\alpha_2,\alpha_3$, DOZZ were able to derive an expression for the residue of $C_\gamma(\alpha_1,\alpha_2,\alpha_3)$ when $s=-n$ with $n$ integer. Let us denote this residue by $\text{Res}_{\sum_{k=1}^3 \alpha_k-2Q= - \gamma n} \: C_\gamma(\alpha_1,\alpha_2,\alpha_3)$. For $s=-n$ given by \eqref{defs} and $\alpha_1,\alpha_2,\alpha_3$ satisfying the bounds \eqref{ThextendedSeibergbounds}, we get (using Fubini to interchange $\E[.]$ and $\int_{\C}$)
\begin{align*}
\text{Res}_{\sum_{k=1}^3 \alpha_k-2Q= - \gamma n} \: C_\gamma(\alpha_1,\alpha_2,\alpha_3) & = 2 \frac{(-\mu)^n}{n!}  \E[\rho(\alpha_1,\alpha_2,\alpha_3)^{n}    ]\\
& = 2 \frac{(-\mu)^n}{n!} \E    [  ( \int_{\C}   \frac{|x|_+^{\gamma(\alpha_1+\alpha_2+\alpha_3)}}{ |x|^{\gamma \alpha_1}  |x-1|^{\gamma \alpha_2}  }  
  M_\gamma(d^2x) )^n ]  \\
&= 2 \frac{(-\mu)^n}{n!}\int_{\C^n}  \prod_{j=1}^n  \frac{|x_j|_+^{\gamma(\alpha_1+\alpha_2+\alpha_3)}}{ |x_j|^{\gamma \alpha_1}  |x_j-1|^{\gamma \alpha_2}}   \E[  \prod_{j=1}^n e^{  \gamma X(x_j)- \frac{\gamma^2}{2}\E[ X(x_j)^2]}  ]   \prod_{j=1}^n \frac{d^2x_j}{|x_j|_+^4}  \\
& = 2 \frac{(-\mu)^n}{n!}\int_{\C^n}  \prod_{j=1}^n  \frac{|x_j|_+^{\gamma(\alpha_1+\alpha_2+\alpha_3)+\gamma^2 (n-1)-4}}{ |x_j|^{\gamma \alpha_1}  |x_j-1|^{\gamma \alpha_2}}  \prod_{i<j} \frac{1}{|x_i-x_j|^{\gamma^2}}  \prod_{j=1}^n d^2x_j  \\
& = 2 \frac{(-\mu)^n}{n!} \int_{\C^n}  \prod_{j=1}^n  \Big ( \frac{1}{ |x_j|^{\gamma \alpha_1}  |x_j-1|^{\gamma \alpha_2}} \Big ) \prod_{i<j} \frac{1}{|x_i-x_j|^{\gamma^2}}  \prod_{j=1}^n d^2x_j,
\end{align*}
where in the last line we have used the fact that $\gamma (\alpha_1+\alpha_2+\alpha_3)= 4- \gamma^2(n-1)$. In conclusion, we have
\begin{equation}\label{linkCFTDF}
\text{Res}_{\sum_{k=1}^3 \alpha_k-2Q= - \gamma n} \: C_\gamma(\alpha_1,\alpha_2,\alpha_3)= 2 \frac{(-\mu)^n}{n!} \int_{\C^n}  \prod_{j=1}^n  \Big ( \frac{1}{ |x_j|^{\gamma \alpha_1}  |x_j-1|^{\gamma \alpha_2}} \Big ) \prod_{i<j} \frac{1}{|x_i-x_j|^{\gamma^2}}  \prod_{j=1}^n d^2x_j
\end{equation}
Notice that the computation which lead to \eqref{linkCFTDF} can be made rigorous by working with the regularized $X_\epsilon$ instead of $X$ and then taking the limit as $\epsilon$ goes to $0$. The key observation of DOZZ is that the last integral above is a famous Dotsenko-Fateev integral and it has an explicit expression in terms of the Gamma function. More specifically, Dotsenko-Fateev \cite{DotFat} found that\footnote{Such integrals were computed by Dotsenko-Fateev in order to give formulas for the three point structure constants of the so-called minimal models of CFT.}
\begin{align*}
&  2 \frac{(-\mu)^n}{n!} \int_{\C^n}  \prod_{j=1}^n  \Big ( \frac{1}{ |x_j|^{\gamma \alpha_1}  |x_j-1|^{\gamma \alpha_2}} \Big ) \prod_{i<j} \frac{1}{|x_i-x_j|^{\gamma^2}}  \prod_{j=1}^n d^2x_j  \\
& =  2 \Big ( \frac{- \pi \mu}{ l(-\frac{\gamma^2}{4}) }  \Big )^n \frac{\prod_{j=1}^n  l(-  \frac{j\gamma^2}{4})}{\prod_{k=0}^{n-1} ( l( \frac{\gamma \alpha_1}{2} +  \frac{k\gamma^2}{4}) l( \frac{\gamma \alpha_2}{2} +  \frac{k\gamma^2}{4}) l( \frac{\gamma \alpha_3}{2} +  \frac{k\gamma^2}{4})  )  }  \\
&:= I_n(\alpha_1,\alpha_2,\alpha_3),   \\
\end{align*}
where recall our convention that $l(x)=\frac{\Gamma (x)}{\Gamma (1-x)}$. Let us set $\bar{\alpha}= \sum_{k=1}^3 \alpha_k$. We therefore have the following relation when $\bar{\alpha}-2Q= - \gamma n$
\begin{align*}
\frac{I_{n-1}(\alpha_1+\gamma ,\alpha_2,\alpha_3)} {I_n(\alpha_1,\alpha_2,\alpha_3)} & = - \frac{l(-\frac{\gamma^2}{4})}{\pi \mu}  \frac{l(\frac{\gamma \alpha_1}{2}) l( \frac{\gamma \alpha_1}{2}+\frac{\gamma^2}{4} )  l(\frac{\gamma}{4}  (2 \alpha_2  + (n-1) \gamma) l(\frac{\gamma}{4}  (2 \alpha_3  + (n-1) \gamma)   )  }{l(- \frac{n \gamma^2}{4})l( \frac{\gamma}{4} (2 \alpha_1+n \gamma) ) }  \\
&=  - \frac{l(-\frac{\gamma^2}{4})}{\pi \mu}  \frac{l(\frac{\gamma \alpha_1}{2}) l( \frac{\gamma \alpha_1}{2}+\frac{\gamma^2}{4} )  l(\frac{\gamma}{4}  (\bar{\alpha}-2 \alpha_1  - \gamma)) }{l( \frac{\gamma}{4}(\bar{\alpha}-2Q))l( \frac{\gamma}{4} (\bar{ \alpha}-2 \alpha_2) ) l( \frac{\gamma}{4} (\bar{ \alpha}-2 \alpha_3) ) },   \\
\end{align*}
where in the last line, we have used the obvious relation $l(x)l(1-x)=1$ and we have made the substitution $n\gamma = 2Q-\bar{\alpha}$. Now, DOZZ argued that by ``analytic continuation" this relation should hold without any restriction on $\alpha_1,\alpha_2,\alpha_3$\footnote{On the level of mathematics, this is not rigorous since for fixed $\alpha_2,\alpha_3$ there is only a finite number of values for $\alpha_1$ such that $-s$ is an integer along with the bounds \eqref{ThextendedSeibergbounds}; one can not analytically continue in a unique way a function defined on a finite set of values}, i.e. one should have
\begin{equation}\label{shiftbyDOZZ}
\frac{C_\gamma(\alpha_1+\gamma ,\alpha_2,\alpha_3)} {C_\gamma(\alpha_1,\alpha_2,\alpha_3)} =  - \frac{l(-\frac{\gamma^2}{4})}{\pi \mu}  \frac{l(\frac{\gamma \alpha_1}{2}) l( \frac{\gamma \alpha_1}{2}+\frac{\gamma^2}{4} )  l(\frac{\gamma}{4}  (\bar{\alpha}-2 \alpha_1  - \gamma)) }{l( \frac{\gamma}{4}(\bar{\alpha}-2Q))l( \frac{\gamma}{4} (\bar{ \alpha}-2 \alpha_2) ) l( \frac{\gamma}{4} (\bar{ \alpha}-2 \alpha_3) ) }
\end{equation}
This enabled them to guess that $C_\gamma(\alpha_1,\alpha_2,\alpha_3)$ could be given by the DOZZ formula since by using \eqref{shiftUpsilon} and after a bit of (tedious!) algebra one can show that $C_\gamma^{{\rm DOZZ}}(\alpha_1,\alpha_2,\alpha_3)$ does indeed satisfy relation  \eqref{shiftbyDOZZ}\footnote{Even by the physicists' standards the derivation lacks rigor. To quote Zamolodchikov-Zamolodchikov \cite{ZaZaarxiv}: 
``It should be stressed that the arguments of this section have nothing to do with a derivation. These are rather some motivations and we consider the expression proposed as a guess which we try to support in the subsequent sections."}.

\section{Lecture 2: The reflection coefficient and the quantum sphere}

In the second lecture, we will explain how to construct the two point correlation functions of LCFT; in the context of LCFT, these correlations are called the reflection coefficient. We will define these correlations denoted $R$ as the following limit for all $\alpha \in (\frac{\gamma}{2},Q)$
\begin{equation}\label{firstdefcorrel2pts}
4 R(\alpha) := \underset{\epsilon \to 0}{\lim} \: \epsilon \:  C_\gamma(\alpha,\epsilon,\alpha)
\end{equation}

In lecture 2, $\alpha$ will always denote a real number in the open interval $(\frac{\gamma}{2},Q)$. We will see  that the reflection coefficient $R(\alpha)$ emerges in the analysis of the tail behaviour of the following random variable
\begin{equation}\label{defIalpha}
I(\alpha):=  \int_{|x| \leq 1}  \frac{1}{|x|^{\alpha \gamma}}   M_\gamma (d^2x).   
\end{equation}

Therefore, we will first explain how to derive tail expansions for $I(\alpha)$. To describe the tail expansion of $I(\alpha)$, we must first recall material introduced in the paper by Duplantier-Miller-Sheffield \cite{DMS}. This is no coincidence since the paper \cite{DMS} develops a theory of random surfaces with two marked points $0$ and $\infty$ hence at the level of LCFT this corresponds to the two point correlation functions. We will show that (the unit volume version of) $R(\gamma)$ is nothing but the partition function of the theory developed in \cite{DMS}. We start by recalling a classical result on drifted Brownian motion by Williams:

\subsection{A few reminders on drifted Brownian motion}

The following decomposition lemma due to Williams (see \cite{Williams})  will be useful in the study of $I(\alpha)$:

\begin{lemma}\label{lemmaWilliams}
Let $(B_s-\nu s)_{s \geq 0}$ be a Brownian motion with negative drift ($\nu >0$) and let $M= \sup_{s \geq 0} (B_s-\nu s)$.  Then conditionally on $M$ the law of the path   $(B_s-\nu s)_{s \geq 0}$ is given by the joining of two  independent paths:

\begin{itemize}
\item
A Brownian motion $((B_s^1+\nu s))_{0 \leq s \leq \tau_M}$ with positive drift $\nu >0$ run until its hitting time $\tau_M$ of $M$. 
\item
$(M+{B}^2_t- \nu t)_{t\geq 0}$ where ${B}^2_t- \nu t$  is a  Brownian motion with negative drift conditioned to stay negative.
\end{itemize}

Moreover, one has the following time reversal property for all $C>0$ (where $\tau_C$ denotes the hitting time of $C$) 
\begin{equation}\label{Timereversal}
(B_{\tau_C-s}^1+\nu (\tau_C-s)-C)_{0 \leq s  \leq \tau_C}\stackrel{law}=  (\tilde{B}_s- \nu s)_{0 \leq s \leq L_{-C}}
\end{equation}
where $(\tilde{B}_s- \nu s)_{s \geq 0}$ is a Brownian motion with drift $-\nu$ conditioned to stay negative and $ L_{-C}$ is the last time  $(\tilde{B}_s- \nu s)$ hits $-C$.  

\end{lemma}

\begin{remark}\label{fundamentalremark}
As a consequence of the above lemma, one can also deduce that the process $(\tilde{B}_{L_{-C}+s}- \nu (L_{-C}+s)+C)_{s \geq 0}$ is equal in distribution to $(\tilde{B}_s- \nu s)_{s \geq 0}$. 
\end{remark}
 
 Usually, Lemma \ref{lemmaWilliams} is stated without the time reversal property; however, this time reversal property is very useful in analysing the tails of $I(\alpha)$.

\subsection{Tail expansion of GMC}

In order to study the tail expansion of $I(\alpha)$, we first need to introduce some notations. We recall basic material introduced in \cite{DMS}.  We first define the process $ \mathcal{B}^\alpha_s$ 
\begin{equation*}
 \mathcal{B}^\alpha_s = \left\{
 \begin{array}{ll}
  B^\alpha_{-s} & \text{if } s < 0\\
    \bar{B}^\alpha_{s} & \text{if } s >0 \end{array} \right.
\end{equation*}
where $B^{\alpha}_s,\bar B^{\alpha}_s$ are  two independent Brownian motions with negative drift $\alpha-Q$ and conditioned to stay negative. Let $Y$ be an independent field with  covariance 
\begin{equation}\label{covlateral}
\E[  Y(s,\theta) Y(t,\theta') ]  = \ln \frac{e^{-s}\vee e^{-t}}{|e^{-s}e^{i \theta} - e^{-t} e^{i \theta'} |}.
\end{equation}
The covariance of the field $Y$ has the same log singularity on the diagonal than the GFF $X$ and belongs to $\mathcal{S}'(\C)$. Following \cite{DMS}, we call the field $Y$ the lateral noise. We also introduce the Gaussian chaos measure with respect to $Y$ 
\begin{equation}\label{chaosmeasure}
 N_\gamma(ds d\theta):= e^{\gamma Y(s,\theta)-\frac{ \gamma^2 }{2}E[Y(s,\theta)^2]}  ds d\theta.
\end{equation}
Just like for the GFF $X$, the chaos measure is defined via a limiting procedure. We also introduce its total mass process
\begin{equation}\label{DefZ}
Z_s = \int_0^{2 \pi} e^{\gamma Y(s,\theta)-\frac{\gamma^2 }{2}E[Y(s,\theta)^2]} d\theta.  
\end{equation}
This is a slight abuse of notation since the process $(Z_s)_{s \in \R}$ is not a function (for $\gamma \geq \sqrt{2}$) but rather a generalized function. With this convention, notice that $Z_s ds$ is stationary i.e. for all $t$ the equality $(Z_{t+s})_{s \in \R}=(Z_s)_{s \in \R}$  holds in distribution. Finally, We set
\begin{equation}\label{rhoalpha}
\rho(\alpha)= \int_{-\infty}^\infty  e^{   \gamma \mathcal{B}_s^\alpha  } Z_s ds.  
\end{equation}
We define the {\it unit volume reflection coefficient} $\bar{R}(\alpha)$  for $\alpha \in (\frac{\gamma}{2},Q)$ by the following formula
\begin{equation}\label{defunitR}
 \bar{R}(\alpha)=  \E[  \rho(\alpha)^{\frac{2}{\gamma}  (Q-\alpha)}].
\end{equation}
 $\bar{R}(\alpha)$  is  indeed well defined as one can show that
 \begin{equation}
\E[  \rho(\alpha)^{p}]< \infty \: \: \Longleftrightarrow \: \: p < \frac{4}{\gamma^2}.
\end{equation} 
and for $\alpha \in (\frac{\gamma}{2},Q)$ we have $\frac{2(Q-\alpha)}{\gamma}<\frac{4}{\gamma^2}$. We can also define the full reflection coefficient by the formula for $\alpha \in (\frac{\gamma}{2},Q)$
 \begin{equation}\label{fullreflection}
 R(\alpha)= \mu^{\frac{2(Q-\alpha)}{\gamma}}   \Gamma (-\frac{2(Q-\alpha)}{\gamma})  \frac{2(Q-\alpha)}{\gamma} \bar R(\alpha)
 \end{equation}
The function $R(\alpha)$ has a divergence at the points $\alpha=\frac{2}{\gamma}- \frac{n-1}{2} \gamma $ with $n \geq 1$ because of the $\Gamma$ function entering the definition. Now, we can state the following theorem on the tail expansion:

\begin{theorem}\label{TailIalpha}
For all $\alpha \in (\frac{\gamma}{2},Q)$ there exists $\eta >0$ (depending on $\alpha$) such that 
\begin{equation*}
\P(    I(\alpha) >t)= \frac{\bar{R}(\alpha)}{t^{\frac{2}{\gamma}(Q-\alpha)  }}  + O(   \frac{1}{t^{\frac{2}{\gamma}(Q-\alpha) +\eta }} )
\end{equation*}

\end{theorem}

\noindent
\emph{Sketch of proof}:

We consider the polar decomposition of $X$ around $0$ as introduced in \cite{DMS}. Consider 
\begin{equation}\label{polardecomp}
X_{e^{-s}}(0):=\frac{1}{2\pi} \int_0^{2\pi} X(e^{-s} e^{i \theta}) d\theta 
\end{equation} 
the circle average around $0$ with radius $e^{-s}$ for $s \geq 0$. A simple computation yields that $B_s:=X_{e^{-s}}(0)$ is a standard Brownian Motion starting from the origin at $s=0$. We have the decomposition
\begin{equation*}
X(e^{-s} e^{i \theta})=X_{e^{-s}}(0)+Y(s,\theta)
\end{equation*}
where $Y$ is the lateral noise with covariance given by \eqref{covlateral}. Hence, one has by the change of variable $x=e^{-s}e^{i \theta}$ 
\begin{align*}
I(\alpha) & =   \int_{|x| \leq 1}  \frac{1}{|x|^{\alpha \gamma}}  e^{\gamma X(x)-\frac{\gamma^2}{2}\E[X(x)^2]} d^2x \\
&  = \int_{0}^\infty \int_0^{2 \pi}  e^{-(2-\alpha)s}   e^{\gamma X(e^{-s} e^{i \theta})-\frac{\gamma^2}{2}\E[X(e^{-s} e^{i \theta})^2]} d\theta ds   \\
&  = \int_{0}^\infty \int_0^{2 \pi}  e^{-(2-\alpha)s}   e^{\gamma B_s -\frac{\gamma^2}{2}s  + \gamma Y(s,\theta) - \frac{\gamma^2}{2} \E[Y(s,\theta)^2]} d\theta ds  \\
&  = \int_{0}^\infty  e^{\gamma (B_s-(Q-\alpha)s)} \left (  \int_0^{2 \pi}   e^{\gamma Y(s,\theta) - \frac{\gamma^2}{2} \E[Y(s,\theta)^2]} d\theta  \right )  ds 
\end{align*}
where in the last line we have used the fact that $Q=\frac{\gamma}{2}+\frac{2}{\gamma}$.
In conclusion we get
\begin{equation}\label{ialphadef}
I(\alpha) =\int_0^\infty  e^{\gamma (B_s-(Q-\alpha)s)}Z_s ds.
\end{equation}
Notice that the computation which lead to \eqref{ialphadef} can be made rigorous by working with the regularized $X_\epsilon$ instead of $X$ and then taking the limit as $\epsilon$ goes to $0$. 
Now, we may apply Lemma \ref{lemmaWilliams} to \eqref{ialphadef}. Let $M=\sup_{s\geq 0}(B_s-(Q-\alpha)s)$ and $L_{-M}$ be the last time $(B^\alpha_s)_{s \geq 0}$ hits $-M$. Then
\begin{align}
&\int_0^\infty e^{\gamma (B_s-(Q-\alpha)s)}Z_s ds\stackrel{law}=e^{\gamma M}\int_{-L_{-M}}^\infty e^{\gamma \mathcal{B}_s^{\alpha}}Z_{s+L_{-M}}ds\stackrel{law}=e^{\gamma M}\int_{-L_{-M}}^\infty e^{\gamma \mathcal{B}_s^{\alpha}}Z_{s} ds\label{will11}
\end{align}
where we used stationarity of the process $Z_s$ (and independence of $Z_s$ and $B_s$). We claim that the tail behaviour of $I(\alpha)$ coincides with that of
\begin{equation*}
J(\alpha)= e^{\gamma M }   \int_{-\infty}^\infty  e^{\gamma \mathcal{B}^\alpha_s}  Z_s ds.
\end{equation*}
The distribution of $M$ is well known (see section 3.5.C in the textbook \cite{KaraSh} for instance):
\begin{equation}\label{tailofmax}
 \P (e^{\gamma M}>t ) = \frac{1}{t^{\frac{2(Q-\alpha)}{\gamma}}} ,  \quad  t\geq 1
\end{equation}
Now the following inequalities hold almost surely
\begin{equation}\label{fund2inequalities}
 e^{\gamma M}\int_{0}^\infty e^{\gamma \mathcal{B}_s^{\alpha}}Z_{s} ds \leq e^{\gamma M}\int_{-L_{-M}}^\infty e^{\gamma \mathcal{B}_s^{\alpha}}Z_{s} ds  \leq e^{\gamma M}\int_{-\infty}^\infty e^{\gamma \mathcal{B}_s^{\alpha}}Z_{s} ds.
\end{equation}
A simple scaling argument using the independence of $(\mathcal{B}_s^{\alpha}, Z_{s})$ with $M$ and the fact that $\int_{-\infty}^\infty e^{\gamma \mathcal{B}_s^{\alpha}}Z_{s} ds$ and $\int_{0}^\infty e^{\gamma \mathcal{B}_s^{\alpha}}Z_{s} ds$ have moments of order $p$ for all $p<\frac{4}{\gamma^2}$  (recall that $\frac{2(Q-\alpha)}{\gamma} <\frac{4}{\gamma^2} $) ensures that the left hand side and the right hand side of \eqref{fund2inequalities} has a tail of the form $\frac{C}{t^{\frac{2(Q-\alpha)}{\gamma}}}$. Indeed, 
\begin{equation}\label{tailofmax1}
\P   ( J(\alpha) >t  ) \underset {t \to \infty}{\sim}   \E  [   \rho(\alpha)^{\frac{2(Q-\alpha)}{\gamma}} ] t^{-\frac{2(Q-\alpha)}{\gamma}}.
\end{equation} 
and similarly we get for the $e^{\gamma M}\int_{0}^\infty e^{\gamma \mathcal{B}_s^{\alpha}}Z_{s} ds$ term
\begin{equation*}
\P   \left (  e^{\gamma M}\int_{0}^\infty e^{\gamma \mathcal{B}_s^{\alpha}}Z_{s} ds   >t  \right ) \underset {t \to \infty}{\sim}   \E  \left [ \left (  \int_{0}^\infty e^{\gamma \mathcal{B}_s^{\alpha}}Z_{s} ds   \right )  ^{\frac{2(Q-\alpha)}{\gamma}} \right ] t^{-\frac{2(Q-\alpha)}{\gamma}}.
\end{equation*} 

This entails that the tail of $e^{\gamma M} \int_{-L_{-M}}^\infty e^{\gamma \mathcal{B}_s^{\alpha}}Z_{s} ds$ is bounded below and above by a constant times $\frac{1}{t^{\frac{2(Q-\alpha)}{\gamma}}}$. On the event that $M$ is bounded, the tail of $e^{\gamma M}\int_{-L_{-M}}^\infty e^{\gamma \mathcal{B}_s^{\alpha}}Z_{s} ds$ is bounded by the tail of $C \int_{-\infty}^\infty e^{\gamma \mathcal{B}_s^{\alpha}}Z_{s} ds$ and hence by $\frac{1}{t^p}$ for all $p<\frac{4}{\gamma^2}$. Therefore, the tail of $e^{\gamma M}\int_{-L_{-M}}^\infty e^{\gamma \mathcal{B}_s^{\alpha}}Z_{s} ds$ is concentrated on the event ``$M$ is big". Therefore, the tail of $e^{\gamma M}\int_{-L_{-M}}^\infty e^{\gamma \mathcal{B}_s^{\alpha}}Z_{s} ds$ is roughly the same as the tail of   
$J(\alpha)= e^{\gamma M }  \rho(\alpha)$. Now, one can conclude by using \eqref{tailofmax1}.

\qed

\subsection{The two point correlation of LCFT and the quantum sphere}

With the previous tail estimates, we can now give a construction of the reflection coefficient of LCFT. Indeed, we have:

\begin{lemma}\label{defR}
For all $\alpha \in (\frac{\gamma}{2},Q)  \setminus \cup_{n\geq 1} \lbrace \frac{2}{\gamma} -\frac{n-1}{2} \gamma \rbrace $, the following limit holds
\begin{equation*}
\underset{\epsilon \to 0}{\lim} \:  \epsilon \: C_\gamma(\alpha,\epsilon,\alpha)= 4 R(\alpha)
\end{equation*}
\end{lemma}

\noindent
\emph{Sketch of proof}:

Recall that 
\begin{equation}\label{expression3pointstructalpha}
C_\gamma(\alpha,\epsilon,\alpha)=
2 \: \mu^{\frac{2(Q-\alpha)}{\gamma}-\frac{\epsilon}{\gamma}}  \gamma^{-1}\Gamma(-\frac{2(Q-\alpha)}{\gamma}+\frac{\epsilon}{\gamma})  \E (\rho(\alpha,\epsilon,\alpha)^{\frac{2(Q-\alpha)}{\gamma}-\frac{\epsilon}{\gamma}})
\end{equation}
where
\begin{equation*}
\rho(\alpha,\epsilon,\alpha)= 
 \int_{\C} 
  \frac{|x|_+^{\gamma(2 \alpha+\epsilon)}}{ |x|^{\gamma \alpha}  |x-1|^{\gamma \epsilon}  }  
  M_\gamma(d^2x).
\end{equation*}
When $\epsilon$ is small we have
\begin{equation*}
\rho(\alpha,\epsilon,\alpha) \approx \int_{\C} 
  \frac{|x|_+^{2 \gamma \alpha}}{ |x|^{\gamma \alpha}   }  
  M_\gamma(d^2x)= I(\alpha)+ I'(\alpha)
\end{equation*}
where
\begin{equation*}
I'(\alpha)= \int_{|x| \geq 1} 
  \frac{|x|_+^{2 \gamma \alpha}}{ |x|^{\gamma \alpha}   }  
  M_\gamma(d^2x)
\end{equation*}
By the change of variable $x  \rightarrow \frac{1}{x}$, $I'(\alpha)$ has the same distribution as $I(\alpha)$. The variables $I(\alpha)$ and $I'(\alpha)$ are ``weakly correlated'' since the tail of $I(\alpha)$ is the same as the tail of $\int_{|x| \leq r} 
  \frac{1}{ |x|^{\gamma \alpha}   }  
  M_\gamma(d^2x)$
for all $r>0$. Therefore, one can show that the tails of $I(\alpha)$ and $I'(\alpha)$ add up, i.e. one can show that there exists $\eta>0$ such that
\begin{equation*}
\P( I(\alpha)+I'(\alpha) >t ) = 2 \frac{\bar{R}(\alpha)}{t^{\frac{2(Q-\alpha)}{\gamma}}}+ O( \frac{1}{t^{\frac{2(Q-\alpha)}{\gamma}+\eta}}  ).
\end{equation*} 
This implies that
\begin{equation*}
\epsilon \: \E[   (I(\alpha)+I'(\alpha)) ^{\frac{2(Q-\alpha)}{\gamma}-\frac{\epsilon}{\gamma}} ] \underset{\epsilon \to 0}{\longrightarrow} 2 \gamma \frac{2(Q-\alpha)}{\gamma} \bar{R}(\alpha)
\end{equation*}
hence the result in view of \eqref{expression3pointstructalpha}.

\qed

Now, we can state an immediate corollary to the DOZZ formula. A straightforward computation based on the shift relations \eqref{shiftUpsilon} yields that
\begin{equation*}
\underset{\epsilon \to 0}{\lim} \:  \epsilon \: C_\gamma^{{\rm DOZZ}}(\alpha,\epsilon,\alpha)= 4 R^{{\rm DOZZ}}(\alpha)
\end{equation*}
 where 
 \begin{equation}\label{Thereflectionrelation}
 R^{{\rm DOZZ}}(\alpha)= - (\pi \: \mu \:  l(\frac{\gamma^2}{4})  )^{\frac{2 (Q -\alpha)}{\gamma}} \frac{\Gamma(-\frac{\gamma (Q-\alpha)}{2})}{\Gamma(\frac{\gamma (Q-\alpha)}{2})} \:  \frac{\Gamma(-\frac{2 (Q-\alpha)}{\gamma })}{\Gamma(\frac{2 (Q-\alpha)}{\gamma})}.
\end{equation}

Therefore, combining Lemma \ref{defR} with Theorem \ref{theoremDOZZ} yields:
\begin{corollary}[Kupiainen, Rhodes, Vargas]\label{theoremRDOZZ}
For all $\alpha \in (\frac{\gamma}{2},Q)$, one has $R(\alpha)= R^{{\rm DOZZ}}(\alpha)$.
\end{corollary}

Finally, let us explain why the (unit volume) reflection coefficient $\bar{R}(\alpha)$ is the partition function of a random measure called the (unit volume) $\alpha$-quantum sphere and introduced by Duplantier-Miller-Sheffield \cite{DMS}. First, we recall the definition of  the $\alpha$-quantum sphere:

\begin{definition}[Duplantier-Miller-Sheffield]
The (unit volume) $\alpha$-quantum sphere is the unit volume random measure $\mu(ds d\theta)$ defined on the cylinder $\R \times [0,2\pi]$  by
\begin{equation}\label{defquantumsphere}
\E[  F(\mu(ds d\theta)) ]=   \frac{\E[    F(  \frac{e^{   \gamma \mathcal{B}_s^\alpha  }  N_\gamma(ds d\theta)}{\rho(\alpha)} )    \rho(\alpha)^{\frac{2}{\gamma}  (Q-\alpha)}      ]}{  \bar{R}(\alpha) }
\end{equation}
\end{definition} 

Therefore, the $\alpha$-quantum sphere is not exactly the measure  $e^{   \gamma \mathcal{B}_s^\alpha  }  N_\gamma(ds d\theta)$ divided by total mass since one must reweight the underlying probability measure by a Radon-Nikodym variable $\rho(\alpha)^{\frac{2}{\gamma}  (Q-\alpha)}$. In fact, the presence of this variable in the definition is a manifestation of the Liouville potential $e^{\gamma \varphi(x)}$ in \eqref{defLiouvilleaction} (with $\mu>0$); recall that we will explain in lecture 3 that the $\E[.]$ term in the correlation function \eqref{Firstdefcorrel} for $N \geq 3$ comes from the Liouville potential and hence the  $ \E[\rho(\alpha)^{\frac{2}{\gamma}  (Q-\alpha)}]$ term in the reflection coefficient comes from the Liouville potential.


\begin{remark}
In fact, Duplantier-Miller-Sheffield \cite{DMS} consider the above random measure in the case $\alpha=\gamma$ in the space of {\bf quantum surfaces}. Essentially, this corresponds to considering the above random measure up to translations along the horizontal or vertical axis of the cylinder. Along the same lines as lemma \ref{defR}, one can show that the unit volume Liouville measure constructed in David-Kupiainen-Rhodes-Vargas \cite{DKRV} with marked points $(\gamma,0), (\epsilon,1), (\gamma,\infty)$ converges when $\epsilon$ goes to $0$ towards  the $\gamma$-quantum sphere in the space of quantum surfaces. For more on the relation between the two approaches, one can also have a look at the paper by Aru-Huang-Sun \cite{aru}.
\end{remark}

\section{Lecture 3: The path integral construction of LCFT}

Now, we explain how to interpret the path integral formulation \eqref{defcorrelations} in order to get definition \eqref{Firstdefcorrel}. We will also explain the KPZ relation \eqref{KPZformula}. Finally, we will state a precise conjecture relating random planar maps to LCFT: this conjecture can be seen as a faithful mathematical formulation of the KPZ conjecture in physics.

\subsection{The path integral construction}
We will consider the metric $g(x)= \frac{1}{|x|_+^4}$; the curvature of the metric is a measure and is given (with generalized function notation) by $R_g(x) g(x) d^2x=-\Delta \ln g(x) d^2x = 4 \nu(d^2x)$ where $\nu$ is  the uniform measure on the circle of center $0$ and radius $1$ (normalized such that $\int_\C\nu(d^2x)= 2 \pi $). Therefore the GFF $X$ is normalized such that it has average $0$ with respect to the curvature $\int_{\C} X(x) R_g(x)g(x)d^2x=0 $. Let us consider 
\begin{equation*}
L^2(\S^2):= \lbrace   \varphi; \; \int_{\C}  \varphi(x) g(x) d^2x < \infty \rbrace
\end{equation*}
the space of standard square integrable functions from  $\S^2$ to $\R$. The standard $H^1(\S^2)$ space is then
\begin{equation*}
H^1(\S^2):=  \lbrace   \varphi; \; \int_{\C}  \varphi(x) g(x) d^2x +   \int_{\C}  |\nabla \varphi(x)|^2  d^2x< \infty  \rbrace 
\end{equation*}  
where $\nabla$ is the standard Gradient in $\C$ (with respect to the flat Euclidean metric). Let $(\varphi_j)_{j \geq 1}$ be the eigenvector basis for $-\Delta_g$, i.e. 
\begin{equation*}
-\frac{1}{g(x)} \Delta \varphi_j(x) = \lambda_{j} \varphi_j(x).
\end{equation*}
normalized to have $L^2(\S^2)$ norm equal to $1$: $\int_{\C}  \varphi_j(x)^2 g(x) d^2x=1$.
Then every function in $\varphi  \in L^2(\S^2)$ can be decomposed in a unique way on the orthonormal basis $(1,(\varphi_j)_{j \geq 1})$
\begin{equation}\label{decompvarphi}
\varphi= c+ \sum_{j \geq 1} c_j \varphi_j 
\end{equation} 
where  for all $j \geq 1$ $c_j=\int_{\C}  f(x) \varphi_j(x) g(x) d^2x $. It is known that the ``Lebesgue measure" $D\varphi$ on $L^2(\S^2)$ does not exist mathematically\footnote{in the physics literature, this formal measure is called the Free Field measure.} but at a formal level, it is natural to write for a function $F$ defined on $L^2(\S^2)$ that
\begin{equation}\label{firstformal}
\int_{L^2(\S^2)}  F( \varphi  ) D\varphi = \int_{\R} \int_{\R^{\mathbb{N}^*}}  F(c + \sum_{j \geq 1} c_j \varphi_j )   \: dc \prod_{j=1}^\infty dc_j
\end{equation}
where $dc$ and each $dc_j$ is the standard Lebesgue measure on $\R$. If $\varphi$ has decomposition \eqref{decompvarphi} then 
\begin{equation*}
\frac{1}{4\pi}\int_{\mathbb{S}^2}  |\nabla_g \varphi(x) |^2 g(x)\, d^2 x= \frac{1}{4\pi }\sum_{j=1}^\infty c_j^2 \lambda_j
\end{equation*} 
hence this leads to the following formal definition
\begin{equation}\label{secondformal}
\int_{L^2(\S^2)} F(\varphi) e^{-\frac{1}{4\pi}\int_{\mathbb{S}^2}  |\nabla_g \varphi(x) |^2 g(x)\, d^2 x }  D\varphi = \int_{\R} \int_{\R^{\mathbb{N}^*}}  F(c + \sum_{j \geq 1} c_j \varphi_j )   \: dc \left ( \prod_{j=1}^\infty e^{  -\frac{c_j^2 \lambda_j}{4\pi} }dc_j  \right )
\end{equation}
Let us stress that the two previous definitions \eqref{firstformal} and \eqref{secondformal} are not meant to be rigorous. However, one can make sense of the previous definition \eqref{secondformal} using probability theory. First, let us make the change of variable $u_j=\frac{c_j \sqrt{\lambda_j}}{\sqrt{2 \pi}}$ in \eqref{secondformal} which leads to (at the formal level)  
\begin{equation*}
 \int_{\R} \int_{\R^{\mathbb{N}^*}}  F(c + \sum_{j \geq 1} c_j \varphi_j )   \: dc \left ( \prod_{j=1}^\infty e^{  -\frac{c_j^2 \lambda_j}{4\pi} }dc_j  \right )= C \int_{\R} \int_{\R^{\mathbb{N}^*}}  F(c + \sqrt{2 \pi}\sum_{j \geq 1} u_j \frac{\varphi_j}{\sqrt{\lambda_j}} )   \: dc \left ( \prod_{j=1}^\infty e^{  -\frac{u_j^2 }{2} }\frac{du_j}{\sqrt{2 \pi}}  \right ) 
\end{equation*}
where the ``constant" $C$ has the following formal definition $C=\prod_{j=1}^\infty (2 \pi (\lambda_j)^{-1/2})$. This constant can be interpreted as $(\text{Det}'(\Delta_g))^{-1/2}$ where $\text{Det}'(\Delta_g)$ is the determinant of the Laplacian (defined mathematically via regularization techniques) but we will disregard it in the sequel. Now, for any i.i.d. sequence $(\epsilon_j)_{j \geq 1}$ of standard centered Gaussian variables the sum $\sqrt{2\pi} \sum_{j \geq 1} \epsilon_j \frac{\varphi_j}{\sqrt{\lambda_j}} $ converges in $\mathcal{S'(\C)}$ to the GFF $X$; hence this leads to the following {\bf rigorous} definition for any function $F$ defined on $\mathcal{S'(\C)}$
\begin{equation}
\int F(\varphi) e^{-\frac{1}{4\pi}\int_{\mathbb{S}^2}  |\nabla_g \varphi(x) |^2 g(x)\, d^2 x }  D\varphi := \int_{\R}  \E[  F(X+c) ]  \: dc.
\end{equation}
Let us stress that since the GFF $X$ is defined in $\mathcal{S'(\C)}$ the underlying space where the above measure lives is actually a space of distributions and not a space of functions and in particular not $L^2(\S^2)$. Since by construction $\int_{\C} X(x)  R_g(x)g(x)d^2x=0$ almost surely and $\int_{\C} R_g(x)g(x)d^2x=8 \pi$, this leads naturally to the following definition 
\begin{equation}
\int F(\varphi) e^{-\frac{1}{4\pi}\int_{\mathbb{S}^2}  |\nabla_g \varphi(x) |^2 g(x)d^2x   -\frac{1}{4\pi}  \int_{\mathbb{S}^2}   R_g(x) \varphi(x) g(x)\, d^2 x }  D\varphi := \int_{\R}  e^{-2Qc} \E[  F(X+c) ]  \: dc
\end{equation}
Let us introduce the {\bf Liouville field} 
\begin{equation}\label{defphi}
\phi= X+\frac{Q}{2} \ln g+c
\end{equation}
and consider the measure $\E^Q$\footnote{The measure has infinite volume and hence is not a probability measure.}
\begin{equation}
\E^Q[  F(\phi) ]= \int_{\R}  e^{-2Qc} \E[  F(X+\frac{Q}{2} \ln g+c) ]  \: dc
\end{equation}

The above measure satisfies the following change of coordinate formula which is classical in the physics literature on LCFT:

\begin{lemma}[change of coordinate formula]\label{lemchange}
For any M\"obius transform, we have for all $F$
\begin{equation}\label{changevar}
\E^Q[  F(\phi \circ \psi+Q \ln |\psi'|) ]= \E^Q[  F(\phi) ]
\end{equation}
\end{lemma}

\begin{remark}
In our construction of LCFT, the change of coordinate formula \eqref{changevar} for $\phi$ is a theorem whereas in \cite{DMS} it is taken as the basis of the definition of a so-called quantum surface. 
\end{remark}

\noindent
\emph{Sketch of proof}:

The proof of \eqref{changevar} relies on the following key identity in distribution
\begin{equation}\label{keyidentity}
X \circ \psi -\frac{1}{2 \pi}\int_{0}^{2 \pi}  (X \circ \psi)(e^{i\theta})  d\theta =  X. 
\end{equation}
Let us set $\mathcal{X}=\frac{1}{2 \pi} \int_{0}^{2 \pi}  (X \circ \psi)(e^{i\theta}) d\theta $; one has using Fubini
\begin{align*}
\E^Q[  F(\phi \circ \psi+Q \ln |\psi'|) ]  & =  \int_{\R}  e^{-2Qc} \E[  F(X\circ \psi+\frac{Q}{2} \ln ( (g\circ \psi) |\psi'|^2)+c) ]  \: dc \\
& =  \E  \left [ \int_{\R}  e^{-2Qc}   F( X\circ \psi+\frac{Q}{2} \ln ((g \circ \psi) |\psi'|^2) +c)  \: dc \right  ] \\
& =  \E \left [ \int_{\R}  e^{-2Qc'} e^{2Q\mathcal{X} -2Q^2 \E[\mathcal{X}^2]}  F( X\circ \psi-\mathcal{X}+Q\E[\mathcal{X}^2]  +\frac{Q}{2} \ln ( (g \circ \psi)|\psi'|^2) +c') ]  \: dc' \right ] \\
& =   \int_{\R}  e^{-2Qc'} \E  [ e^{2Q\mathcal{X} -2Q^2 \E[\mathcal{X}^2] }  F( X\circ \psi-\mathcal{X}+Q\E[\mathcal{X}^2] +\frac{Q}{2} \ln ((g \circ \psi)|\psi'|^2) +c') ] \:  dc', \\
\end{align*}
where we made the change of variable $c'=c+\mathcal{X}-Q \E[\mathcal{X}^2]$. By using the Girsanov Theorem \ref{th:Girsanov} (in the appendix), we get for all $c'$
\begin{align*}
 &  \E  [ e^{2Q\mathcal{X} -2Q^2 \E[\mathcal{X}^2] }  F( ((X\circ \psi)(x)-\mathcal{X}+Q\E[\mathcal{X}^2] +\frac{Q}{2} \ln ((g \circ \psi)(x)|\psi'(x)|^2) +c')_{x \in \C}) ] \\
 & =  \E  [   F( (X\circ \psi)(x)+ v(x)+c')_{x \in \C} ]
\end{align*} 
where $v(x)= 2Q \E[  (X\circ \psi)(x) \mathcal{X}    ]-Q\E[\mathcal{X}^2]+ \frac{Q}{2} \ln ((g \circ \psi)(x)|\psi'(x)|^2) $. Now, one can conclude by a (lengthy!) computation that 
\begin{equation*}
v(x)=\frac{Q}{2} \ln g(x)
\end{equation*}
This yields the result.
\qed

Let $\phi_\epsilon$ be the circle average approximation of $\phi$, namely $\phi_\epsilon(x)= \frac{1}{2\pi} \int_0^{2\pi}  \phi(x+\epsilon e^{i \theta}) d\theta  $. We have $\phi_\epsilon=X_\epsilon+\frac{Q}{2}(\ln g)_\epsilon+c$ where $(\ln g)_\epsilon$ denotes the circle average of $\ln g$. We consider the associated ``Vertex operator" 
\begin{equation*}
V_{\alpha,\epsilon}(x)= \epsilon^{\frac{\alpha^2}{2}} e^{\alpha\phi_\epsilon(x)}
\end{equation*}
We have for $x \in \C$ that $\E[ X_\epsilon(x)^2 ]= \ln \frac{1}{\epsilon}-\frac{1}{2} \ln g(x)+o(1)$ ($o(1)$ is with respect to $\epsilon$) hence
\begin{equation*}
V_{\alpha,\epsilon}(x)= e^{\alpha c} e^{\alpha X_\epsilon(x)- \frac{\alpha^2}{2}  \E[ X_\epsilon(x) ]} g(x)^{ \frac{\alpha Q}{2}-\frac{\alpha^2}{4} }(1+o(1)) 
\end{equation*}

For $\alpha=\gamma$, one gets $\frac{\gamma Q}{2}-\frac{\gamma^2}{4}=1$ and therefore
\begin{equation*}
V_{\gamma,\epsilon}(x)= e^{\gamma c} \frac{e^{\gamma X_\epsilon(x)- \frac{\gamma^2}{2}  \E[ X_\epsilon(x) ]}}{|x|_+^4}(1+o(1)) 
\end{equation*}
hence we get the following convergence in the space of Radon measures
\begin{equation*}
V_{\gamma,\epsilon}(x) d^2x \underset{\epsilon \to 0}{\rightarrow} e^{\gamma c} M_{\gamma} (d^2x).
\end{equation*}

In view of the previous considerations, it is then natural to define the correlations by the following formula
\begin{equation}\label{defcorrelpathintegral}
\langle \prod_{i=1}^N V_{\alpha_i} (z_i) \rangle_{\gamma,\mu}  := \underset{\epsilon \to 0}{\lim} \: \langle \prod_{i=1}^N V_{\alpha_i,\epsilon} (z_i) \rangle_{\gamma,\mu}
\end{equation}
where 
\begin{equation*}
\langle \prod_{i=1}^N V_{\alpha_i,\epsilon} (z_i) \rangle_{\gamma,\mu}  := 2�\: \E^Q[ e^{-\mu \int_{\C} V_{\gamma,\epsilon}(x) d^2x}  \prod_{k=1}^N V_{\alpha_k,\epsilon} (z_k)   ]
\end{equation*}

Let us gather this construction in a proposition:
\begin{proposition}
The limit in \eqref{defcorrelpathintegral} exists and is equal to the definition \eqref{Firstdefcorrel} introduced in lecture 1.
\end{proposition}

\proof

First, we get by Fubini (interchanging $\E[.]$ and $dc$)
\begin{align*}
& \E^Q[ e^{-\mu \int_{\C} V_{\gamma,\epsilon}(x) d^2x}  \prod_{k=1}^N V_{\alpha_k,\epsilon} (z_k)   ]  \\
& = \E[   \int_{\R}    e^{-2 Q c}  \prod_{k=1}^N  \epsilon^{\frac{\alpha_k^2}{2}} e^{  \alpha_k  (X_\epsilon(z_k)+ \frac{Q}{2} (\ln g)_\epsilon(z_k) + c)  }   e^{-\mu e^{\gamma c} \epsilon^{\frac{\gamma^2}{2}} \int_{\C} e^{  \gamma  (X_\epsilon(x)+ \frac{Q}{2} (\ln g)_\epsilon(x) )  }  d^2x}  dc     ]  \\
& = \gamma^{-1}\int_0^\infty  u^{s-1}e^{-\mu u} du \:  \E \left [      \frac{ \prod_{k=1}^N \epsilon^{\frac{\alpha_k^2}{2}} e^{  \alpha_k  (X_\epsilon(z_k)+ \frac{Q}{2} (\ln g)_\epsilon (z_k))  }}  { \Big ( \epsilon^{\frac{\gamma^2}{2}}  \int_{\C} e^{  \gamma  (X_\epsilon(x)+ \frac{Q}{2} (\ln g)_\epsilon (x) )  }  d^2x \Big)^s }       \right ], \\
\end{align*}
where $s=\frac{\sum_{k=1}^N \alpha_k-2Q}{\gamma}$ and using the change of variable $u= \mu e^{\gamma c} \epsilon^{\frac{\gamma^2}{2}} \int_{\C} e^{  \gamma  (X_\epsilon(x)+ \frac{Q}{2} (\ln g)_\epsilon (x)  }  d^2x$ in the last line. Now, using the Girsanov Theorem \ref{th:Girsanov} (in the appendix) applied to $\mathcal{X}= \sum_{k=1}^N \alpha_k  X_\epsilon(z_k)$, this yields
\begin{align*}
& 2 \: \E^Q[ e^{-\mu \int_{\C} V_{\gamma,\epsilon}(x) d^2x}  \prod_{k=1}^N V_{\alpha_k,\epsilon} (z_k)   ]  \\
& =2 \mu^{-s} \gamma^{-1}\Gamma (s)  \prod_{k=1}^N  e^{   (\frac{\alpha_k Q}{2}-\frac{\alpha_k^2}{4}) (\ln g)_\epsilon(z_k)  } \prod_{i<j} e^{   \alpha_i \alpha_j \E[X_\epsilon( z_i)X_\epsilon(z_j)]   } \\
&  \times \E \left [      \Big ( \epsilon^{\frac{\gamma^2}{2}}  \int_{\C} e^{  \gamma  (X_\epsilon(x)+\sum_{k=1}^N \alpha_k\E[ X_\epsilon(x)X_\epsilon(z_k) ] + \frac{Q}{2} (\ln g)_\epsilon (x) )   }  d^2x \Big)^{-s}        \right ] \\
& \underset{\epsilon \to 0}{\rightarrow}  2 \mu^{-s} \gamma^{-1}\Gamma(s)
 \prod_{i < j} \frac{1}{|z_i-z_j|^{\alpha_i \alpha_j}}\E \left [  \left (  \int_{\C}  F(x,{\bf z}) M_\gamma(d^2x)  \right )^{-s}  \right ]   \\
\end{align*}
where we have used the fact that $\E[X_\epsilon( x)X_\epsilon(y)] \underset{\epsilon \to 0}{\rightarrow} \ln \frac{1}{|x-y|}+ \ln |x|_++\ln |y|_+$ (which is equal to $\ln \frac{1}{|x-y|}-\frac{1}{4} \ln g(x)-\frac{1}{4} \ln g(y) $) and recall that $F(x,{\bf z})=\prod_{k=1}^N \left ( \frac{ |x|_+}{|x-z_k|}  \right )^{\gamma \alpha_k}$. 

\qed

We can now prove the KPZ formula:

\begin{proposition}[KPZ formula]
\begin{equation}\label{KPZformulaprop}
\langle \prod_{k=1}^N V_{\alpha_k}(\psi(z_k))    \rangle=  \prod_{k=1}^N |\psi'(z_k)|^{-2 \Delta_{\alpha_k}}     \langle \prod_{k=1}^N V_{\alpha_k}(z_k)     \rangle
\end{equation}  
where  $\Delta_{\alpha}=\frac{\alpha}{2}(Q-\frac{\alpha}{2})$ is called the conformal weight.

\end{proposition}

\proof

Let us fix $\epsilon>0$. By using the change of coordinate lemma \ref{lemchange} with
\begin{equation*}
F(\phi)=  e^{-\mu \int_{\C} V_{\gamma,\epsilon}(x) d^2x}  \prod_{k=1}^N V_{\alpha_k,\epsilon} (z_k)
\end{equation*}
we get
\begin{equation}\label{exactidentity}
\E^Q[ e^{-\mu \epsilon^{\frac{\gamma^2}{2}}\int_{\C}e^{ \gamma (\phi \circ \psi)_\epsilon(x)+ Q \ln |\psi'(x)|_\epsilon  } d^2x}  \prod_{k=1}^N \epsilon^{\frac{\alpha_k^2}{2}}  e^{ \alpha_k (\phi \circ \psi)_\epsilon(z_k)+ Q \ln |\psi'(z_k)|_\epsilon  }  ] = \E^Q[ e^{-\mu \int_{\C} V_{\gamma,\epsilon}(x) d^2x}  \prod_{k=1}^N V_{\alpha_k,\epsilon} (z_k)   ]
\end{equation}
where $\ln |\psi'(x)|_\epsilon$ is the circle average approximation to $\ln |\psi'(x)|$. When $\epsilon$ is small, we have  $\ln |\psi'(x)|_\epsilon\approx \ln |\psi'(x)|$ hence relation \eqref{exactidentity} leads to
\begin{align}
& \Big ( \prod_{k=1}^N |\psi'(z_k)|^{2 \Delta_{\alpha_k}} \Big )\E^Q[ e^{-\mu \int_{\C} (|\psi'(x)|\epsilon)^{\frac{\gamma^2}{2}}e^{ \gamma (\phi \circ \psi)_\epsilon(x) } |\psi'(x)|^{2\Delta_{\gamma}}d^2x}  \prod_{k=1}^N (|\psi'(z_k)|\epsilon)^{\frac{\alpha_k^2}{2}}  e^{ \alpha_k (\phi \circ \psi)_\epsilon(z_k)  }  ]  \nonumber  \\
& \approx \E^Q[ e^{-\mu \int_{\C} V_{\gamma,\epsilon}(z) d^2z}  \prod_{k=1}^N V_{\alpha_k,\epsilon} (z_k)   ],   \label{approxpsi}
\end{align}
 Now, we have by definition of the circle average
\begin{equation*}
(\phi \circ \psi)_\epsilon(x) \approx \phi_{|\psi'(x)|\epsilon}(\psi(x)). 
\end{equation*}
In particular, using the fact that $\Delta_{\gamma}=1$ this yields that 
\begin{align*}
& \int_{\C} (|\psi'(x)|\epsilon)^{\frac{\gamma^2}{2}}e^{ \gamma (\phi \circ \psi)_\epsilon(x) } |\psi'(x)|^{2\Delta_\gamma}d^2x \\
& = \int_{\C} (|\psi'(x)|\epsilon)^{\frac{\gamma^2}{2}}e^{ \gamma (\phi \circ \psi)_\epsilon(x) } |\psi'(x)|^2d^2x \\
& \approx \int_{\C} (|\psi'(x)|\epsilon)^{\frac{\gamma^2}{2}}e^{ \gamma \phi_{|\psi'(x)|\epsilon}( \psi (x)) } |\psi'(x)|^2d^2x \\
&  = \int_{\C} (|\psi'(\psi^{-1}(u))|\epsilon)^{\frac{\gamma^2}{2}}e^{ \gamma \phi_{|\psi'(\psi^{-1}(u))|\epsilon}( u) } d^2u\\
& \underset{\epsilon \to 0}{\rightarrow} M_\gamma(\C)\\ 
\end{align*}
by a change of variable $ \epsilon_u=  |\psi'(\psi^{-1}(u))|\epsilon$ as $\epsilon$ goes to $0$. Then, using $(\phi \circ \psi)_\epsilon(z_k) \approx \phi_{|\psi'(z_k)|\epsilon}(\psi(z_k)) $ and a change of variable $\epsilon_k=  |\psi'(z_k)|\epsilon$ at each vertex operator, we get that
\begin{align*}
& \Big ( \prod_{k=1}^N |\psi'(z_k)|^{2 \Delta_{\alpha_k}} \Big )\E^Q[ e^{-\mu \int_{\C} (|\psi'(x)|\epsilon)^{\frac{\gamma^2}{2}}e^{ \gamma (\phi \circ \psi)_\epsilon(x) } |\psi'(x)|^2d^2x}  \prod_{k=1}^N (|\psi'(z_k)|\epsilon)^{\frac{\alpha_k^2}{2}}  e^{ \alpha_k (\phi \circ \psi)_\epsilon(z_k)  }  ]   \\
&  \underset{\epsilon \to 0}{\rightarrow} \frac{1}{2}( \prod_{k=1}^N |\psi'(z_k)|^{2 \Delta_{\alpha_k}} \Big ) \langle \prod_{k=1}^N V_{\alpha_k}(\psi(z_k))    \rangle  \\
\end{align*}
hence the conclusion in view of \eqref{approxpsi}.
\qed

\subsection{The KPZ conjecture and equation}

Here we explain the so-called KPZ conjecture (after Knizhnik-Polyakov-Zamolodchikov \cite{cf:KPZ}) relating the scaling limit of random fields on a regular lattice and the scaling limit of random fields on a random lattice, i.e. random planar maps (equipped with a conformal structure). The KPZ conjecture we state in this subsection is not to be confused with other KPZ relations \cite{Ar,BGRV,cf:DuSh,Rnew10} which have appeared in the probability literature and which relate the dimensions of a set with respect to the Euclidian metric to the dimension of the same set with respect to the measure $M_\gamma(d^2x)$. One should see these KPZ relations based on dimensions as a geometrical interpretation of the KPZ conjecture we state below.

Consider a random field $\sigma(x) d^2x$\footnote{We adopt the notation of generalized functions.} living in $\mathcal{S}'(\C)$ which is expected to describe the scaling limit of an observable of a critical statistical physics model considered on a regular lattice. The scaling limit is to be considered as the mesh of the lattice goes to $0$ and we suppose that it corresponds to a CFT with central charge $c_M < 1$. According to CFT, the correlations of the field (if it is ``primary") will behave like a conformal tensor: for all M\"obius transforms $\psi$
\begin{equation*}
\langle \prod_{k=1}^N \sigma (\psi(z_k))    \rangle=  \prod_{k=1}^N |\psi'(z_k)|^{-2 \Delta_{\sigma}}     \langle \prod_{k=1}^N \sigma(z_k)     \rangle
\end{equation*}  
where $\Delta_\sigma$ is a real number called the conformal weight. In a 1988 seminal paper, Knizhnik-Polyakov-Zamolodchikov (KPZ) \cite{cf:KPZ} argued that the scaling limit $ \Sigma(x) d^2x$ of the same observable on a random planar map conformally embedded in the sphere should factorize into $\sigma$ and an independent part depending on LCFT\footnote{The probability of a planar map is chosen to be proportional to the partition function (at critical temperature) of the model of statistical physics living on the map.}. More precisely, suppose that $\Sigma(x)d^2x $ is considered in the sphere by conformally mapping $3$ points at random on the map to three fixed points of the complex plane $z_1,z_2,z_3 \in \C$. We denote $\Sigma^{z_1,z_2,z_3}(x) d^2x$ the scaling limit to stress that it depends on $z_1,z_2,z_3$ . KPZ argued that one should roughly have $\Sigma^{z_1,z_2,z_3}(x)= \sigma(x) V_\alpha(x)$ for some $\alpha$ where $V_\alpha$ is a vertex operator of LCFT with parameter $\gamma$ ($\alpha$ and $\gamma$ are different except when $\Delta_\sigma=0$). One can determine $\gamma$ by solving the following equation discovered by Polyakov in his 1981 seminal paper \cite{Pol}
\begin{equation}\label{centralcharge0}
c_M+c_L=26
\end{equation}
where $c_M$ is the central charge describing $\sigma$  and $c_L=1+6Q^2$ is the central charge of LCFT. The equation \eqref{centralcharge0} has a unique solution $\gamma \in (0,2)$ if $c_M<1$: this solution is $\gamma= \frac{\sqrt{25-c_M}-\sqrt{1-c_M}}{\sqrt{6}}$.  In terms of correlations, we get the following explicit mathematical conjecture: the field $\Sigma^{z_1,z_2,z_3}(x) d^2x$ has the following correlation structure
\begin{equation}\label{defcorrelsKPZ}
\langle \prod_{k=1}^N \Sigma^{z_1,z_2,z_3}(x_k)     \rangle:=  \langle \prod_{k=1}^N \sigma(x_k)     \rangle  \frac{\langle V_\gamma(z_1) V_\gamma(z_2) V_\gamma (z_3) \prod_{k=1}^N V_\alpha(x_k)     \rangle_{\gamma,\mu}}{
  \langle V_\gamma(z_1) V_\gamma(z_2) V_\gamma (z_3)     \rangle_{\gamma,\mu}}   \quad \text{ ({\bf KPZ conjecture})}
\end{equation}  
where the $V_\gamma (z_k)$ terms in the above relation correspond to the embedding of the map. By construction, the field $\Sigma^{z_1,z_2,z_3}(x) d^2x$ is conformally invariant, i.e. one must have for all $A \subset \S^2$
\begin{equation*}
\int_{\psi(A)}  \Sigma^{\psi(z_1),\psi(z_2),\psi(z_3)}(u) d^2u=  \int_{A}  \Sigma^{z_1,z_2,z_3}(x) d^2x
\end{equation*}
This conformal invariance property enforces that
\begin{equation*}
\langle \prod_{k=1}^N \Sigma^{\psi(z_1),\psi(z_2),\psi(z_3)}(\psi(x_k))     \rangle=   \prod_{k=1}^N |\psi'(z_k)|^{-2 }     \langle \prod_{k=1}^N   \Sigma^{z_1,z_2,z_3}(x_k)     \rangle
\end{equation*}  
Applying this to the right hand side of \eqref{defcorrelsKPZ} leads to the quadratic {\bf KPZ equation}:
\begin{equation}\label{KPZquadratic}
 \Delta_{\sigma}+\Delta_{\alpha}=1 \quad \quad \quad \quad \quad \quad \quad \quad \text{ ({\bf KPZ equation})}
\end{equation} 
where recall that $\Delta_{\alpha}= \frac{\alpha}{2}(Q-\frac{\alpha}{2})$. If $\Delta_{\sigma}>1-\frac{Q^2}{4} $ there is a unique $\alpha<Q$ which solves \eqref{KPZquadratic}.

\subsubsection*{The particular case $\sigma=1$}
In the special case where $\sigma=1$, the field  $\Sigma^{z_1,z_2,z_3}(x) d^2x$ is the conjectured scaling limit of the volume form of the random planar map conformally embedded in the sphere. By conformal invariance, we can suppose that $z_1=0, z_2=1, z_3=\infty$. In this special case, one has  $\Delta_{\sigma}=0$ and solving the KPZ equation \eqref{KPZquadratic} leads to $\alpha=\gamma$. In fact, we can integrate  \eqref{defcorrelsKPZ} on an Euclidean ball $B$ to give an explicit expression for the law of the volume form $\Sigma^{0,1,\infty}(x) d^2x$. Let us set introduce the following notations (which generalize \eqref{defrho3})
\begin{equation*}
\rho(\alpha_1,\alpha_2,\alpha_3) [d^2x]=  \frac{|x|_+^{\gamma(\alpha_1+\alpha_2+\alpha_3)}}{ |x|^{\gamma \alpha_1}  |x-1|^{\gamma \alpha_2}  }  
  M_\gamma(d^2x),   \quad  \quad  \rho(\alpha_1,\alpha_2,\alpha_3):=\rho(\alpha_1,\alpha_2,\alpha_3)[\C] 
\end{equation*}
After some reverse engineering on the Girsanov formula, we get the following expression
\begin{equation*}
\int_{B^n}  \langle V_\gamma(0) V_\gamma(1) V_\gamma (\infty) \prod_{k=1}^N V_\gamma(x_k)     \rangle_{\gamma,\mu}  \prod_{k=1}^N d^2x_k= \frac{2}{\gamma} \mu^{-\frac{(N+3)\gamma-2Q}{\gamma} }  \Gamma \left ( \frac{(N+3)\gamma-2Q}{\gamma} \right )  \E[  ( \frac{\rho(\gamma,\gamma,\gamma) [B] }{ \rho(\gamma,\gamma,\gamma) } )^N      \rho(\gamma,\gamma,\gamma)^{- \frac{3 \gamma-2Q}{\gamma}}  ]  
\end{equation*}      
Hence, one gets the following definition for the conjectured scaling limit of the random planar map volume form: 
\begin{definition} [{\bf Liouville volume form}]
The Liouville volume form (or measure) $ \nu_{L}(d^2x)=\Sigma^{0,1,\infty}(x) d^2x$ is a random measure defined by 
\begin{equation*}
\E[   F( \nu_{L} (d^2x) )]= \E \left [  F \left ( \xi \frac{\rho(\gamma,\gamma,\gamma) [d^2x] }{ \rho(\gamma,\gamma,\gamma) }   \right )      \rho(\gamma,\gamma,\gamma)^{- \frac{3 \gamma-2Q}{\gamma}}  \right ]/ \E[ \rho(\gamma,\gamma,\gamma)^{- \frac{3 \gamma-2Q}{\gamma}}  ]
\end{equation*} 
where $\xi$ is an independent variable with a Gamma density $\frac{\mu^{\frac{3 \gamma-2Q}{\gamma}}}{\Gamma (  \frac{3 \gamma-2Q}{\gamma}    )}t^{ \frac{3 \gamma-2Q}{\gamma}-1}e^{- \mu t} dt$.

\end{definition}

\section{Lecture 4: sketch of proof of the DOZZ formula}

In this lecture, we explain the main ideas behind the proof of the DOZZ formula stated as Theorem \ref{theoremDOZZ} in lecture 1. We denote $\mathcal{A}_N$ the (convex) set of real numbers $(\alpha_1,\cdots, \alpha_N)$ satisfying condition \eqref{ThextendedSeibergbounds}.

In order to keep these notes rather concise and introductory, we will make the following assumptions:

\begin{itemize}
\item
{\bf  Analycity}: for all $z_1, \cdots, z_N \in \C$ the map 
\begin{equation*}
(\alpha_1, \cdots, \alpha_N) \mapsto \frac{\langle \prod_{k=1}^N V_{\alpha_k} (z_k) \rangle_{\gamma,\mu}}{\Gamma(\frac{\sum_{k=1}^N\alpha_k-2Q}{\gamma})}
\end{equation*} 
is analytic on the set $\mathcal{A}_N$. 

\item
{\bf  Reflection relation}: The probabilistically defined reflection coefficient $R$ given by \eqref{fullreflection} (see lecture 2) satisfies $R(\alpha)=R^{{\rm DOZZ}}(\alpha)$ for $\alpha \in (\frac{\gamma}{2},Q)$. In particular, this implies that $R$ can be extended to a meromorphic function on $\C$ satisfying two remarkable shift equations
\begin{equation}\label{firstshift}
R(\alpha)= - \mu \pi \frac{R(\alpha+\frac{\gamma}{2})}{ l(-\frac{\gamma^2}{4}) l(\frac{\gamma\alpha}{2})  l(2+\frac{\gamma^2}{4}- \frac{\gamma \alpha}{2})}
\end{equation} 
and
\begin{equation}\label{secondshift}
R(\alpha)=   -  \tilde{\mu} \pi \frac{R(\alpha+\frac{2}{\gamma})}{l(-\frac{4}{\gamma^2})l(\frac{2 \alpha}{\gamma}) l(2+\frac{4}{\gamma^2}-\frac{2 \alpha}{\gamma})} 
\end{equation} 
where $ \tilde{\mu}= \frac{(\mu \pi l(\frac{\gamma^2}{4})  )^{\frac{4}{\gamma^2} }}{\pi l(\frac{4}{\gamma^2})}$ is the dual cosmological constant. One can notice the duality $\frac{\gamma}{2} \leftrightarrow \frac{2}{\gamma}, \; \; \mu \leftrightarrow \tilde{\mu}= \frac{(\mu \pi l(\frac{\gamma^2}{4})  )^{\frac{4}{\gamma^2} }}{\pi l(\frac{4}{\gamma^2})}$  in the above relation.
\end{itemize}

In fact, proving both assumptions is non trivial, especially the {\bf  Reflection relation}. Now, we proceed with a sketch of the proof. The proof is based on two main ingredients: the BPZ equations along with crossing symmetry considerations and rigorous Operator Product expansions (OPE).

\subsection{Sketch of the proof using the BPZ differential equations and crossing symmetry}
In this subsection, we explain how to use the BPZ differential equations for degenerate field insertions in 4 point correlation functions to prove the DOZZ formula. The proof will also be based on rigorous OPE that we will explain in the next subsection. 
The four point correlation function is fixed by the M\"obius invariance \eqref{KPZformula} up to a single function depending on the cross ratio of the points. For later purpose we label the points from $0$ to $3$ and consider the weights $\alpha_1, \alpha_2, \alpha_3$ fixed:
\begin{align}
 \nonumber \langle    
   \prod_{k=0}^3 V_{\alpha_k}(z_k)  \rangle_{\gamma,\mu} 
 & = |z_3-z_0|^{- 4 \Delta_{0}}  |z_2-z_1|^{ 2 (\Delta_3-\Delta_2-\Delta_1-\Delta_0)}|z_3-z_1|^{2(\Delta_2+\Delta_0 -\Delta_3 -\Delta_1  )} \\ &\times |z_3-z_2|^{2 (\Delta_1+\Delta_0-\Delta_3-\Delta_2)} G_{\alpha_0}\left ( \frac{(z_0-z_1)(z_2-z_3)}{ (z_0-z_3) (z_2-z_1)}  \right ) . \label{confinv}
\end{align}
where $\Delta_k= \Delta_{\alpha_k}$. We can recover $G_{\alpha_0}$ as the following limit
\begin{equation}
G_{\alpha_0}(z)= \lim_{ z_3 \to\infty }|z_3|^{4 \Delta_3} \langle    V_{\alpha_0}  (z)  V_{\alpha_1}(0)  V_{\alpha_2}(1) V_{\alpha_3}(z_3)  \rangle_{\gamma,\mu}   \label{Glimit0}.
\end{equation}
Combining with  \eqref{Firstdefcorrel}  we get
 \begin{equation}
  G_{\alpha_0}(z)=|z|^{-\alpha_0 \alpha_1}   |z-1|^{-\alpha_0 \alpha_2}  \mathcal{T}_{\alpha_0}(z) 
   \label{Glimit1}
\end{equation}
where, setting $s=\frac{\alpha_0+\alpha_1+\alpha_2+\alpha_3-2Q}{\gamma}$, $ \mathcal{T}_{\alpha_0}
  (z)$ is given by 
\begin{align}\label{Tdefi0}
  \mathcal{T}_{\alpha_0}(z) &=
  2 \mu^{-s}  \gamma^{-1}\Gamma(s)  \E [R_{\alpha_0}(z)^{-s}]
 \end{align}
and
\begin{equation}
R_{\alpha_0}(z)=  
\int_{\C} 
\frac{|x|_+^{\gamma \sum_{k=0}^3 \alpha_k}}{ |x-z|^{\gamma\alpha_0}|x|^{\gamma \alpha_1}  |x-1|^{\gamma \alpha_2}  }  M_\gamma(d^2x).
\label{Rdefi0}
\end{equation}
There are two special values of $\alpha_0$ for which the reduced  four point function $ \mathcal{T}_{\alpha_0}(z)$ satisfies a second order differential equation. That such equations are expected in Conformal Field Theory goes back to BPZ \cite{BPZ}. In the  case of LCFT it was proved in \cite{KRV}  that, under suitable assumptions on $\alpha_1,\alpha_2,\alpha_3$, if $\alpha_0\in\{-\frac{\gamma}{2},-\frac{2}{\gamma}\}$  then $\mathcal{T}_{\alpha_0}$ is a   solution of a PDE version of the Gauss hypergeometric equation
\begin{equation}\label{hypergeo}
\partial_{z}^2 \caT_{\alpha_0}(z)+  \frac{({ c}-z({ a}+{ b}+1))}{z(1-z)}\partial_z \caT_{\alpha_0}(z) -\frac{{ a}{ b} }{z(1-z)}\caT_{\alpha_0}(z)=0
\end{equation}
where ${ a},{ b},{ c}$ are given by
\begin{align}\label{defabc}
{ a}&= \frac{\alpha_0}{2} (Q-2\alpha_0-\alpha_1-\alpha_2-\alpha_3)-\frac{1}{2},  \quad { b}=\frac{\alpha_0}{2} (Q-\alpha_1-\alpha_2+\alpha_3)+\frac{1}{2},  \quad c=1+\alpha_0 (Q-\alpha_1).
\end{align}
This equation has two holomorphic solutions defined on  $\mathbb{C} \setminus \lbrace (-\infty,0) \cup (1,\infty) \rbrace$:
 \begin{equation}\label{Fpmdef}
F_{-}(z)= {}_2F_1({ a},{  b},{ c},z), \quad F_{+}(z)= z^{1-{ c}} {}_2F_1(1+{ a}-{ c},1+{ b}-{ c},2-{ c},z)
\end{equation}
where $_2F_1(a,b,c,z)$ is given by the standard hypergeometric series (which can be extended holomorphically on $\mathbb{C} \setminus  (1,\infty) $). 

It is quite remarkable that the space of real valued solutions to \eqref{hypergeo} in $\mathbb{C} \setminus \{0,1\}$ is a 1d space which one can describe in great detail (Lemma 4.4 in \cite{KRV}). More specifically, all solutions are determined up to a multiplicative constant $\lambda\in\R$ as
\begin{equation}\label{Tsolution}
\caT_{\alpha_0}(z)= 
\lambda( | F_{-}(z) |^2+A_\gamma(\alpha_0,\alpha_1 , \alpha_2, \alpha_3) | F_{+}(z) |^2)
\end{equation}
where  the coefficient $A_\gamma(\alpha_0,\alpha_1 , \alpha_2, \alpha_3)$ is given by
\begin{equation}\label{Fundrelation}
A_\gamma(\alpha_0,\alpha_1 , \alpha_2, \alpha_3)=- \frac{\Gamma(c)^2  \Gamma(1-a)  \Gamma(1-b)  \Gamma(a-c+1)  \Gamma(b-c+1) }{ \Gamma(2-c)^2  \Gamma(c-a)  \Gamma(c-b)  \Gamma (a)  \Gamma(b) }
\end{equation}
provided $c \in \R \setminus \Z$ and $c-a-b \in  \R \setminus \Z$.

Furthermore, the constant $\lambda$ is found by using  the expressions \eqref{expression3pointstruct} and \eqref{Tdefi0} (note that $s$ has a different meaning in these two expressions):
\begin{equation}\label{Fundrelation1}\lambda=\caT_{\alpha_0}(0)=C_\gamma(\alpha_1+\alpha_0,\alpha_2,\alpha_3).
\end{equation}
Hence for $\alpha_0\in\{-\frac{\gamma}{2},-\frac{2}{\gamma}\}$  $\mathcal{T}_{\alpha_0}$ is completely determined in terms of $C_\gamma(\alpha_1+\alpha_0,\alpha_2,\alpha_3)$. 

Now, the main idea behind the proof of the DOZZ formula is to perform an asymptotic expansion around $z=0$ to give another expression  of the coefficient in front of $| F_{+}(z) |^2$ in terms of $C_\gamma(\alpha_1-\alpha_0,\alpha_2,\alpha_3)$.  This will yield a functional relation between $C_\gamma(\alpha_1-\alpha_0,\alpha_2,\alpha_3)$ and $C_\gamma(\alpha_1+\alpha_0,\alpha_2,\alpha_3)$ thanks to \eqref{Tsolution}.

\subsubsection{The case $\alpha_0=-\frac{\gamma}{2}$ and $\alpha_1+\frac{\gamma}{2}<Q$}
If $\sum_{k=1}^3\alpha_k>2Q+\frac{\gamma}{2}$ and $\alpha_1+\frac{\gamma}{2}<Q$ then one can show (see next subsection on OPE) the following expansion around $z=0$
\begin{equation}\label{Fundrelation300}
\caT_{-\frac{\gamma}{2}}(z)= C_\gamma(\alpha_1-\frac{\gamma}{2},\alpha_2,\alpha_3)+ B(\alpha_1)C_\gamma(\alpha_1+\frac{\gamma}{2},\alpha_2,\alpha_3) | z |^{2(1-c)}+o( | z |^{2(1-c)})
\end{equation}
where recall that for $\alpha_0=-\frac{\gamma}{2}$ one has $2(1-c)= \gamma(Q-\alpha_1)$ and $B$ is given by 
\begin{equation}\label{Fundrelation2}
B(\alpha)=-\mu  \frac{\pi}{  l(-\frac{\gamma^2}{4}) l(\frac{\gamma \alpha}{2})  l(2+\frac{\gamma^2}{4}- \frac{\gamma \alpha}{2}) } .
\end{equation}
Hence, in view of \eqref{Tsolution} and \eqref{Fundrelation1}, relations \eqref{Fundrelation300}, \eqref{Fundrelation2} lead to   
\begin{equation}\label{Fundrelation30}B(\alpha_1)C_\gamma(\alpha_1+\frac{\gamma}{2},\alpha_2,\alpha_3)=A_\gamma(-\frac{\gamma}{2}, \alpha_1 , \alpha_2, \alpha_3)C_\gamma(\alpha_1-\frac{\gamma}{2},\alpha_2,\alpha_3)
\end{equation}
 which yields 
the following relation (after some algebra!)
\begin{equation}\label{3pointconstanteqintro}
\frac{C_\gamma(\alpha_1+\frac{\gamma}{2},\alpha_2,\alpha_3)}{C_\gamma(\alpha_1-\frac{\gamma}{2},\alpha_2,\alpha_3)}  
= - \frac{1}{\pi \mu}\frac{l(-\frac{\gamma^2}{4})  l(\frac{\gamma \alpha_1}{2}) l(\frac{\alpha_1\gamma}{2}  -\frac{\gamma^2}{4})  l(\frac{\gamma}{4} (\bar{\alpha}-2\alpha_1- \frac{\gamma}{2}) )   }{l( \frac{\gamma}{4} (\bar{\alpha}-\frac{\gamma}{2} - 2Q)  ) l( \frac{\gamma}{4} (\bar{\alpha}-2\alpha_3-\frac{\gamma}{2} ))  l( \frac{\gamma}{4} (\bar{\alpha}-2\alpha_2-\frac{\gamma}{2} )) }, 
\end{equation}
where $\bar{\alpha}=\alpha_1+\alpha_2+\alpha_3$. Thanks to this last relation, one can analyticaly continue $\alpha \mapsto C_\gamma(\alpha,\alpha_2,\alpha_3)$ to a strip of the form $\R \times [-\eta,\eta]$ with $\eta>0$ (the analytic continuation is  meromorphic with poles). We will also denote this analytic continuation $C_\gamma(\alpha,\alpha_2,\alpha_3)$. 

\subsubsection{the case $\alpha_0=-\frac{\gamma}{2}$ and $\alpha_1+\frac{\gamma}{2}>Q$}
If $\sum_{k=1}^3\alpha_k>2Q+\frac{\gamma}{2}$ and $\alpha_1+\frac{\gamma}{2}>Q$ then one can show (see next subsection on OPE) the following expansion around $z=0$
\begin{equation}\label{Fundrelation3000}
\caT_{-\frac{\gamma}{2}}(z)= C_\gamma(\alpha_1-\frac{\gamma}{2},\alpha_2,\alpha_3)+ R(\alpha_1)C_\gamma(2Q-\alpha_1-\frac{\gamma}{2},\alpha_2,\alpha_3) | z |^{2(1-c)}+o( | z |^{2(1-c)})
\end{equation}
where $R$ is defined by the probabilistic expression \eqref{fullreflection}. Now, recall that we admit that \eqref{firstshift} holds; hence the above expansion can be written
\begin{equation}\label{Fundrelation3000bis}
\caT_{-\frac{\gamma}{2}}(z)= C_\gamma(\alpha_1-\frac{\gamma}{2},\alpha_2,\alpha_3)+B(\alpha_1) R(\alpha_1+\frac{\gamma}{2}) C_\gamma(2Q-\alpha_1-\frac{\gamma}{2},\alpha_2,\alpha_3) | z |^{2(1-c)}+o( | z |^{2(1-c)})
\end{equation}
where $B$ is defined by \eqref{Fundrelation2}. Therefore we get by combining \eqref{Tsolution} and \eqref{Fundrelation1} with \eqref{Fundrelation3000bis}   
\begin{equation}\label{3pointconstanteqintrobis}
\frac{R(\alpha_1+\frac{\gamma}{2})C_\gamma(2Q-\alpha_1-\frac{\gamma}{2},\alpha_2,\alpha_3)}{C_\gamma(\alpha_1-\frac{\gamma}{2},\alpha_2,\alpha_3)}  
= - \frac{1}{\pi \mu}\frac{l(-\frac{\gamma^2}{4})  l(\frac{\gamma \alpha_1}{2}) l(\frac{\alpha_1\gamma}{2}  -\frac{\gamma^2}{4})  l(\frac{\gamma}{4} (\bar{\alpha}-2\alpha_1- \frac{\gamma}{2}) )   }{l( \frac{\gamma}{4} (\bar{\alpha}-\frac{\gamma}{2} - 2Q)  ) l( \frac{\gamma}{4} (\bar{\alpha}-2\alpha_3-\frac{\gamma}{2} ))  l( \frac{\gamma}{4} (\bar{\alpha}-2\alpha_2-\frac{\gamma}{2} )) }. 
\end{equation}
Thanks to this last relation, we see that $\alpha \mapsto R(\alpha)C_\gamma(2Q-\alpha,\alpha_2,\alpha_3)$ provides an analytic continuation of  $C_\gamma(\alpha,\alpha_2,\alpha_3)$ beyond $\alpha=Q$ hence by unicity of analytic continuation one has $R(\alpha)C_\gamma(2Q-\alpha,\alpha_2,\alpha_3)=C_\gamma(\alpha,\alpha_2,\alpha_3)$ (where recall that we denote $C_\gamma(\alpha,\alpha_2,\alpha_3)$ the analytic continuation provided by \eqref{3pointconstanteqintro}).

\subsubsection{The case $\alpha_0=-\frac{2}{\gamma}$}
If $\sum_{k=1}^3\alpha_k>2Q+\frac{2}{\gamma}$  then one can show (see next subsection on OPE) the following expansion around $z=0$
\begin{equation}\label{Fundrelation30000}
\caT_{-\frac{2}{\gamma}}(z)= C_\gamma(\alpha_1-\frac{2}{\gamma},\alpha_2,\alpha_3)+ R(\alpha_1)C_\gamma(2Q-\alpha_1-\frac{2}{\gamma},\alpha_2,\alpha_3) | z |^{2(1-c)}+o( | z |^{2(1-c)})
\end{equation}
where $R$ is defined by the probabilistic expression \eqref{fullreflection} and for $\alpha_0=-\frac{2}{\gamma}$ one has $2(1-c)=\frac{4}{\gamma}(Q-\alpha_1)$. Since we admit \eqref{secondshift}, this is equivalent to  
\begin{equation}\label{Fundrelation30000}
\caT_{-\frac{2}{\gamma}}(z)= C_\gamma(\alpha_1-\frac{2}{\gamma},\alpha_2,\alpha_3)+ \tilde{B}(\alpha_1)R(\alpha_1+\frac{2}{\gamma})C_\gamma(2Q-\alpha_1-\frac{2}{\gamma},\alpha_2,\alpha_3) | z |^{2(1-c)}+o( | z |^{2(1-c)})
\end{equation}
where
\begin{equation}\label{Fundrelation20}
\tilde{B}(\alpha)=-\tilde{\mu}  \frac{\pi}{  l(-\frac{4}{\gamma^2}) l(\frac{2 \alpha}{\gamma})  l(2+\frac{4}{\gamma^2}- \frac{2\alpha}{\gamma}) } .
\end{equation}
with $\tilde{\mu}= \frac{(\mu \pi l(\frac{\gamma^2}{2})  )^{\frac{4}{\gamma^2} }}{ \pi l(\frac{4}{\gamma^2})}$ the dual cosmological constant.
Hence, since $R(\alpha_1+\frac{2}{\gamma})C_\gamma(2Q-\alpha_1-\frac{2}{\gamma},\alpha_2,\alpha_3)=C_\gamma(\alpha_1+\frac{2}{\gamma},\alpha_2,\alpha_3) $, this leads to the following relation (after some algebra!) using \eqref{Tsolution} and \eqref{Fundrelation1}
\begin{equation}\label{3pointconstanteqintrodual}
\frac{C_\gamma(\alpha_1+\frac{2}{\gamma},\alpha_2,\alpha_3)}{C_\gamma(\alpha_1-\frac{2}{\gamma},\alpha_2,\alpha_3)}  
= - \frac{1}{\pi \tilde{\mu}}\frac{l(-\frac{4}{\gamma^2})  l(\frac{2 \alpha_1}{\gamma}) l(\frac{2\alpha_1}{\gamma}  -\frac{4}{\gamma^2})  l(\frac{1}{\gamma} (\bar{\alpha}-2\alpha_1- \frac{2}{\gamma}) )   }{l( \frac{1}{\gamma} (\bar{\alpha}-\frac{2}{\gamma} - 2Q)  ) l( \frac{1}{\gamma} (\bar{\alpha}-2\alpha_3-\frac{2}{\gamma} ))  l( \frac{1}{\gamma} (\bar{\alpha}-2\alpha_2-\frac{2}{\gamma} )) }. 
\end{equation}

\subsubsection{Unicity of the DOZZ shift equations}
We suppose that $\gamma^2 \not \in \Q$; the general case can be deduced from this case by a continuity argument. The DOZZ formula $C_\gamma^{{\rm DOZZ}}(\alpha_1,\alpha_2,\alpha_3)$ satisfies also the two shift equations \eqref{3pointconstanteqintro} and \eqref{3pointconstanteqintrodual}. If we fix $\alpha_2,\alpha_3$ then the function $\varphi(\alpha)=\frac{C_\gamma(\alpha,\alpha_2,\alpha_3)}{C_\gamma^{{\rm DOZZ}}(\alpha,\alpha_2,\alpha_3)}$ has two periods $\gamma$ and $\frac{4}{\gamma}$ hence it is a constant $c_{\alpha_2,\alpha_3}$ (since $\gamma^2 \not \in \Q$). By a similar argument applied to the $\alpha_2$ variable one can show that this constant depends only on $\alpha_3$. Then, working on the $\alpha_3$ variable one can show that the constant $c:= c_{\alpha_2,\alpha_3}$ depends on no variables. This proves the DOZZ formula since the constant $c$ is $1$; indeed, for $\alpha \in (\frac{\gamma}{2},Q)$ 
\begin{equation*}
c= \underset{\epsilon \to 0}{\lim} \: \frac{C_\gamma(\alpha,\epsilon,\alpha)}{C_\gamma^{{\rm DOZZ}}(\alpha,\epsilon,\alpha)}= \frac{R(\alpha)}{R^{{\rm DOZZ}}(\alpha)}=1.
\end{equation*}

\subsection{Operator product expansions}

Finally, we explain how to derive the asymptotic expansions of the previous subsection.

\subsubsection{The case $\alpha_0=-\frac{\gamma}{2}$ and $\alpha_1+\frac{\gamma}{2}<Q$}
We give a sketch of the proof of \eqref{Fundrelation300} in the case 
\begin{equation}\label{extracondition}
\gamma+\frac{1}{\gamma}< \alpha_1+\frac{\gamma}{2}<Q.
\end{equation}
 The case $\alpha_1+\frac{\gamma}{2} \leq \gamma+\frac{1}{\gamma}$ is more involved and we will not discuss this case in these notes. We set $p:=\frac{1}{\gamma}(\sum_{k=1}^3 \alpha_k-\frac{\gamma}{2}-2Q)$ and 
\begin{equation}\label{defineK}
K(z,x)=\frac{|x-z|^{\frac{\gamma^2}{2}}|x|_+^{\gamma (\sum_{k=1}^3\alpha_k-\frac{\gamma}{2})}}{|x|^{\gamma \alpha_1}|x-1|^{\gamma \alpha_2}}
\end{equation}
With these notations, we have
\begin{equation*}
\mathcal{R}_{-\frac{\gamma}{2}}(z)= \int_{\C} 
  K(z,x) M_\gamma(d^2x).
\end{equation*}
We perform the following approximation as $z$ goes to $0$
\begin{equation*}
\mathcal{R}_{-\frac{\gamma}{2}}(z)^{-p}-\mathcal{R}_{-\frac{\gamma}{2}}(0)^{-p} \approx  -p \left ( \mathcal{R}_{-\frac{\gamma}{2}}(z)-\mathcal{R}_{-\frac{\gamma}{2}}(0) \right) \mathcal{R}_{-\frac{\gamma}{2}}(0)^{-p-1}
\end{equation*}
This leads to
\begin{align*}
& \E[ \mathcal{R}_{-\frac{\gamma}{2}}(z)^{-p}]-E[\mathcal{R}_{-\frac{\gamma}{2}}(0)^{-p}]   \\
& \approx  -p \E \left [ \left ( \mathcal{R}_{-\frac{\gamma}{2}}(z)-\mathcal{R}_{-\frac{\gamma}{2}}(0) \right) \mathcal{R}_{-\frac{\gamma}{2}}(0)^{-p-1} \right ]    \\
& = - p \int_\C (K(z,u)-K(0,u))  \E[     \left (    \int_{\C}  K(0,x)e^{\gamma^2 \E[X(x)X(u)]} M_\gamma(d^2x)  \right )^{-p-1}  ]  \frac{du}{|u|_+^4}
\end{align*}
where in the last line we have used Corollary \ref{coro:Girsanov} of the appendix. Now for all continuous function $f$ the following equivalent holds
\begin{equation}\label{lemmaong}
\int_\C (K(z,u)-K(0,u)) f(u)  \frac{du}{|u|_+^4} \underset{z \to 0}{\sim}     |z|^{\gamma (Q-\alpha_1)}  f(0) \frac{\pi}{  l(-\frac{\gamma^2}{4}) l(\frac{\gamma \alpha_1}{2})  l(2+\frac{\gamma^2}{4}- \frac{\gamma \alpha_1}{2}) } 
\end{equation}
Indeed the mass of the above integral concentrates around $u=0$ (when $z$ goes to $0$) and we have 
\begin{align*}
& \int_\C (K(z,u)-K(0,u)) f(u)  \frac{du}{|u|_+^4}   \\
& =  \int_{\C} \frac{(|u-z|^{\frac{\gamma^2}{2}}-|u|^{\frac{\gamma^2}{2}} )|u|_+^{\gamma (\sum_{k=1}^3\alpha_k-\frac{\gamma}{2})}}{|u|^{\gamma \alpha_1}|u-1|^{\gamma \alpha_2}} f(u)  \frac{du}{|u|_+^4}          \\
& \approx  \int_{|u| \leq \frac{1}{2}} \frac{(|u-z|^{\frac{\gamma^2}{2}}-|u|^{\frac{\gamma^2}{2}} )|u|_+^{\gamma (\sum_{k=1}^3\alpha_k-\frac{\gamma}{2})}}{|u|^{\gamma \alpha_1}|u-1|^{\gamma \alpha_2}} f(u)  \frac{du}{|u|_+^4}          \\
&  = |z|^{\gamma (Q-\alpha_1)} \int_{|y| \leq \frac{1}{2|z|}} \frac{(|y-\frac{z}{|z|}|^{\frac{\gamma^2}{2}}-|y|^{\frac{\gamma^2}{2}} )|zy|_+^{\gamma (\sum_{k=1}^3\alpha_k-\frac{\gamma}{2})}}{|y|^{\gamma \alpha_1}|zy-1|^{\gamma \alpha_2}} f(zy)  \frac{dy}{|zy|_+^4}          \\
&  \underset{z \to 0}{\sim} |z|^{\gamma (Q-\alpha_1)} f(0) \int_{\C} \frac{ |y-1|^{\frac{\gamma^2}{2}}-|y|^{\frac{\gamma^2}{2}}  }{|y|^{\gamma \alpha_1} }dy         \\
\end{align*}
 Notice that condition \eqref{extracondition} ensures that $\int_{\C} \frac{ |y-1|^{\frac{\gamma^2}{2}}-|y|^{\frac{\gamma^2}{2}}  }{|y|^{\gamma \alpha_1} }dy$ is well defined. This leads to \eqref{lemmaong} since $\int_{\C} \frac{ |y-1|^{\frac{\gamma^2}{2}}-|y|^{\frac{\gamma^2}{2}}  }{|y|^{\gamma \alpha_1} }dy=\frac{\pi}{  l(-\frac{\gamma^2}{4}) l(\frac{\gamma \alpha_1}{2})  l(2+\frac{\gamma^2}{4}- \frac{\gamma \alpha_1}{2}) } $. Applying \eqref{lemmaong} with $ f(u)= \E \left [\left(    \int_{\C}  K(z,x)e^{\gamma^2 \E[X(x)X(u)] } M_{\gamma}(d^2x)  \right ) ^{-p-1} \right ] $, we get  
\begin{align*}
& \E[ \mathcal{R}_{-\frac{\gamma}{2}}(z)^{-p}]-\E[\mathcal{R}_{-\frac{\gamma}{2}}(0)^{-p}]   \\
& \underset{z \to 0} {\sim}  -p |z|^{\gamma (Q-\alpha_1)}  \frac{\pi}{  l(-\frac{\gamma^2}{4}) l(\frac{\gamma \alpha_1}{2})  l(2+\frac{\gamma^2}{4}- \frac{\gamma \alpha_1}{2}) } \E  \left [   \left (    \int_{\C}  K(0,x)e^{\gamma^2 \E[X(x) X(0)]} M_\gamma(d^2x)  \right )^{-p-1} \right ]. 
\end{align*}
Notice that this implies that
\begin{align*}
& \mathcal{T}_{-\frac{\gamma}{2}}(z)-  \mathcal{T}_{-\frac{\gamma}{2}}(0)  \\
& \underset{z \to 0} {\sim}  - 2  \mu^{-p} \gamma^{-1} p \Gamma(p)  |z|^{\gamma (Q-\alpha_1)}  \frac{\pi}{  l(-\frac{\gamma^2}{4}) l(\frac{\gamma \alpha_1}{2})  l(2+\frac{\gamma^2}{4}- \frac{\gamma \alpha_1}{2}) } \E  \left [   \left (    \int_{\C}  K(0,x)e^{\gamma^2 \E[X(x) X(0)]} M_\gamma(d^2x)  \right )^{-p-1} \right ]. 
\end{align*}

This gives the desired result since $p \Gamma(p)=\Gamma(p+1)$ and 
\begin{equation*}
C_\gamma(\alpha_1+\frac{\gamma}{2},\alpha_2,\alpha_3)= 2  \mu^{-p-1} \gamma^{-1}  \Gamma(p+1)   \E  \left [   \left (    \int_{\C}  K(0,x)e^{\gamma^2 \E[X(x) X(0)]} M_\gamma(d^2x)  \right )^{-p-1} \right ].
\end{equation*}

\subsubsection{The case $\alpha_0=-\frac{\gamma}{2}$ and $\alpha_1+\frac{\gamma}{2}>Q$}
We give a sketch of the proof of \eqref{Fundrelation3000}. We stick to the same notations as the previous proof. The main idea is that the increment 
\begin{equation*}
\E[ \mathcal{R}_{-\frac{\gamma}{2}}(z)^{-p} ]-\E[ \mathcal{R}_{-\frac{\gamma}{2}}(0)^{-p}   ]
\end{equation*}
is ruled by the large tail of the variable   $ \int_{|x| \leq |z|} K(z,x) M_\gamma(d^2x)$. Let $X_{|z|}(0)$ denote the circle average of $X$ on the circle of center $0$ and radius $|z|$ (it has variance $\ln \frac{1}{|z|}$). For $|x| \leq |z|$, the process  $\tilde{X}_z(x)=X(|z|x)-X_{|z|}(0) $ is distributed like $(X(x))_{|x| \leq 1}$.
We have
\begin{align*}
& \int_{|x| \leq |z|} K(z,x) M_\gamma(d^2x)  \\
& = \int_{|x| \leq |z|} \frac{|x-z|^{\frac{\gamma^2}{2}}}{|x|^{\gamma \alpha_1}|1-x|^{\gamma \alpha_2}} e^{\gamma X(x) -\frac{\gamma^2}{2}\E[X(x)^2]  } d^2x  \\
& \approx |z|^{\frac{\gamma^2}{2}}  \int_{|x| \leq |z|} \frac{1}{|x|^{\gamma \alpha_1}} e^{\gamma X(x) -\frac{\gamma^2}{2}\E[X(x)^2]  } d^2x  \\
& =  |z|^{\frac{\gamma^2}{2}} e^{\gamma X_{|z|}(0)+\frac{\gamma^2}{2} \ln |z| } \int_{|x| \leq |z|} \frac{1}{|x|^{\gamma \alpha_1}} e^{\gamma (X(x)-X_{|z|}(0)) -\frac{\gamma^2}{2}\E[(X(x)-X_{|z|}(0))^2]  } d^2x  \\
& =  |z|^{\frac{\gamma^2}{2}-\gamma \alpha_1+2} e^{\gamma X_{|z|}(0)+\frac{\gamma^2}{2} \ln |z| } \int_{|y| \leq 1} \frac{1}{|y|^{\gamma \alpha_1}} e^{\gamma \tilde{X}_z(y) -\frac{\gamma^2}{2}\E[\tilde{X}_z(y)^2]  } d^2x  \\
\end{align*}
We use the notations of lecture 2 in what follows. We denote $\tilde{I}_z(\alpha_1)=\int_{|y| \leq 1} \frac{1}{|y|^{\gamma \alpha_1}} e^{\gamma \tilde{X}_z(y) -\frac{\gamma^2}{2}\E[\tilde{X}_z(y)^2]  } d^2x$ which has same law as $I(\alpha_1)$ given by \eqref{ialphadef}. In the sequel, we will write the following approximation $\tilde{I}_z(\alpha_1) \approx e^{\gamma M} \rho(\alpha_1)$. Recall that $e^{\gamma M}$ has density $\frac{2(Q-\alpha_1)}{\gamma} \frac{1}{v^{\frac{2(Q-\alpha_1)}{\gamma}+1}}$. When $z$ gets small, the variable  $\tilde{I}_z(\alpha_1)$ is roughly independent from $(X(x))_{ |x| \geq |z|}$ hence we get by averaging with respect to $e^{\gamma M}$ and setting $A_z=\int_{|x| \geq |z|} K(z,x) M_\gamma(d^2x)$
\begin{align*}
& \E[ \mathcal{R}_{-\frac{\gamma}{2}}(z)^{-p} ]-\E[ \mathcal{R}_{-\frac{\gamma}{2}}(0)^{-p}   ]  \\
& \approx  \E \left [ \left ( A_z+ \int_{|x| \leq |z|} K(z,x) M_\gamma(d^2x)  \right )^{-p}  -A_z^{-p}    \right  ]  \\
& \approx \E \left [ \left ( A_z+  |z|^{\frac{\gamma^2}{2}-\gamma \alpha_1+2} e^{\gamma X_{|z|}(0)+\frac{\gamma^2}{2} \ln |z|} e^{\gamma M} \rho(\alpha_1) \right )^{-p} -A_z^{-p}    \right  ]  \\
& =  \frac{2(Q-\alpha_1)}{\gamma}  \E \left [  \int_{0}^\infty    \left (    ( A_z+  |z|^{\frac{\gamma^2}{2}-\gamma \alpha_1+2} e^{\gamma X_{|z|}(0)+\frac{\gamma^2}{2}\ln |z|} v \rho(\alpha_1) )^{-p} -A_z^{-p}  \right )   \frac{dv}{v^{\frac{2(Q-\alpha_1)}{\gamma}+1}} \right ]  \\
& = \frac{2(Q-\alpha_1)}{\gamma} \bar{R}(\alpha_1) \left ( \int_{0}^\infty   (  (1+u)^{-p}-1 ) \frac{du}{  u^{\frac{2(Q-\alpha_1)}{\gamma}+1}} \right ) |z|^{\gamma (Q-\alpha_1)}    \E[  A_z^{-p -\frac{2(Q-\alpha_1)}{\gamma} }  e^{2(Q-\alpha_1) X_{|z|}(0)  +2(Q-\alpha_1)^2  \ln |z| }   ]
\end{align*}
where we used the change of variable $u=\frac{ |z|^{\frac{\gamma^2}{2}-\gamma \alpha_1+2} e^{\gamma X_{|z|}+\frac{\gamma^2}{2} \ln |z|} v \rho(\alpha_1)}{A_z}$. By the Girsanov theorem \ref{th:Girsanov}, we get 
\begin{align*}
 & \E[  A_z^{-p -\frac{2(Q-\alpha_1)}{\gamma} }  e^{2(Q-\alpha_1) X_{|z|}(0)  +2(Q-\alpha_1)^2  \ln |z| }   ]   \\
&=  \E \left [  \left ( \int_{|x| \geq |z|} K(z,x) M_\gamma(d^2x) \right ) ^{-p -\frac{2(Q-\alpha_1)}{\gamma} }  e^{2(Q-\alpha_1) X_{|z|}(0)  +2(Q-\alpha_1)^2  \ln |z| }  \right ]   \\
& =  \E \left [  \left ( \int_{|x| \geq |z|} K(z,x)e^{  2 \gamma (Q-\alpha_1) \E[  X(x) X_{|z|}(0)  ] } M_\gamma(d^2x) \right ) ^{-p -\frac{2(Q-\alpha_1)}{\gamma} }  \right  ]   \\ 
& \underset{z \to 0}{\rightarrow} \E \left [  \left ( \int_\C K(0,x)e^{  2 \gamma (Q-\alpha_1) \E[X(x) X(0)] } M_\gamma(d^2x) \right ) ^{-p -\frac{2(Q-\alpha_1)}{\gamma} }   \right ]   \\ 
\end{align*}
hence the result since 
\begin{equation*}
 \int_0^{\infty}    \left (  (1+u )^{-p}-1  \right )   \frac{du}{u^{ \frac{2(Q-\alpha_1)}{\gamma}+1  }}=   \frac{\Gamma(-\frac{2(Q-\alpha_1)}{\gamma})  \Gamma (p+\frac{2(Q-\alpha_1)}{\gamma})}{   \Gamma (p) }
\end{equation*}
and
\begin{equation*}
C_\gamma(2Q-\alpha_1-\frac{\gamma}{2},\alpha_2,\alpha_3)= 2  \mu^{-p-\frac{2(Q-\alpha_1)}{\gamma}} \gamma^{-1}  \Gamma(p+\frac{2(Q-\alpha_1)}{\gamma})   \E  \left [   \left (    \int_{\C}  K(0,x)e^{2\gamma (Q-\alpha_1)  \E[X(x) X(0)]} M_\gamma(d^2x)  \right )^{-p-\frac{2(Q-\alpha_1)}{\gamma}} \right ].
\end{equation*}
\subsubsection{The case $\alpha_0=-\frac{2}{\gamma}$}
This case is similar to the case $\alpha_0=-\frac{\gamma}{2}$ and $\alpha_1+\frac{\gamma}{2}>Q$.
   
  \section{Appendix}

  \subsection{The Girsanov theorem}
  We state the classical Girsanov theorem:
  
 \begin{theorem}{\bf Girsanov theorem }\label{th:Girsanov}

Let $\mathcal{X}$ be some some Gaussian variable which is measurable with respect to the GFF $X$. Let $F$ be some bounded continuous function. Then we have the following identity
\begin{equation*}
\E[  e^{\mathcal{X}-\frac{\E[\mathcal{X}^2]}{2}}   F( (X(x))_{x \in \C}   )   ]= \E[    F( (X(x)   +E[X(x) \mathcal{X}])_{x \in \C}   )   ]
\end{equation*}
\end{theorem}  

  In particular we get the following corollary:
  
\begin{corollary}\label{coro:Girsanov}
 Let $F$ be some bounded continuous function. Then we have the following identity
\begin{equation}\label{coroGirs}
\E[   \left ( \int_{\C}  f(u) M_\gamma(d^2u)   \right )  F( (X(x))_{x \in \C}   )   ]= \int_{\C} f(u)\E[    F( (X(x)   +E[X(x) X(u)])_{x \in \C}   )   ]  \frac{d^2u}{|u|_+^4}
\end{equation}
\end{corollary}  
 
 \proof
 For $\epsilon>0$, if $X_\epsilon$ denotes the circle average of $X$ then by Fubini (interchanging $\E[.]$ and $\int_{\C}$)
\begin{align*}
& \E[   \left ( \int_{\C}  f(u) e^{\gamma X_\epsilon(u)-\frac{\gamma^2}{2} \E[X_\epsilon(u)^2]}\frac{d^2u}{|u|_+^4}   \right )  F( (X(x))_{x \in \C}   )   ]  \\
&  = \int_{\C}  f(u) \E[   e^{\gamma X_\epsilon(u)-\frac{\gamma^2}{2} \E[X_\epsilon(u)^2]}   F( (X(x))_{x \in \C}   )   ]  \frac{d^2u}{|u|_+^4} \\
& = \int_{\C} f(u)  \E[    F( (X(x)   +E[X(x) X_\epsilon(u)])_{x \in \C}   )   ]   \frac{d^2u}{|u|_+^4} \\
\end{align*} 
where in the last line we have used the Girsanov Theorem \ref{th:Girsanov}. In conclusion, we have
\begin{equation*}
 \E[   \left ( \int_{\C}  f(u) e^{\gamma X_\epsilon(u)-\frac{\gamma^2}{2} \E[X_\epsilon(u)^2]}\frac{d^2u}{|u|_+^4}   \right )  F( (X(x))_{x \in \C}   )   ]  � = \int_{\C} f(u)  \E[    F( (X(x)   +E[X(x) X_\epsilon(u)])_{x \in \C}   )   ]   \frac{d^2u}{|u|_+^4} 
\end{equation*} 
We conclude by letting $\epsilon$ go to $0$. 
 \qed


\hspace{10 cm}


\begin{thebibliography}{20}



\bibitem{AGT}
L. F. Alday, D. Gaiotto, and Y. Tachikawa. Liouville Correlation Functions from Four Dimensional Gauge Theories, \emph{Lett. Math. Phys.} \textbf{91}, 167-197 (2010).

\bibitem{Ar}
Aru J.: KPZ relation does not hold for the level lines and $\text{SLE}_\kappa$ flow lines of the Gaussian free field, \emph{Probability Theory and Related Fields} {\bf 163} 3-4, 465-526 (2015).  
\bibitem{aru}
Aru J., Huang Y., Sun X.:  Two perspectives of the 2D unit area quantum sphere and their equivalence,  	arXiv:1512.06190.


\bibitem{BPZ} Belavin A.A., Polyakov A.M., Zamolodchikov A.B. : Infinite conformal symmetry in two-dimensional quantum field theory, \emph{Nuclear Physics B} {\bf 241} (2), 333-380 (1984). 


\bibitem{Ber}
Berestycki N.: An elementary approach to Gaussian multiplicative chaos, \emph{Electronic communications in Probability} {\bf 27}, 1-12 (2017).  

\bibitem{BGRV}
Berestycki N., Garban C., Rhodes R., Vargas V.: KPZ formula derived from Liouville heat kernel, \emph{Journal of the London Mathematical Society} {\bf 94} (1), 186-208 (2016).

 
 

\bibitem{DKRV}
David F., Kupiainen A., Rhodes R., Vargas V.: Liouville Quantum Gravity on the Riemann sphere, \emph{Communications in Mathematical Physics} {\bf 342} (3), 869-907 (2016).

\bibitem{DoOt}
Dorn H., Otto H.-J.: Two and three point functions in Liouville theory, \emph{Nuclear Physics B} {\bf 429} (2), 375-388 (1994).     


\bibitem{DotFat}
Dotsenko V., Fateev V.: Four Point Correlation Functions and the Operator Algebra in the Two-Dimensional Conformal Invariant Theories with the Central Charge $c < 1$, \emph{Nuclear Physics B} {\bf 251}, 691-734 (1985).

\bibitem{Dub0}
Dub\'edat J.: SLE and the Free Field: partition functions and couplings, \emph{Journal of the AMS} \textbf{22} (4), 995-1054 (2009). 


\bibitem{DMS}
Duplantier B., Miller J., Sheffield: Liouville quantum gravity as mating of trees,  	\href{http://arxiv.org/abs/1409.7055}{arXiv:1409.7055}.

\bibitem{cf:DuSh} Duplantier, B., Sheffield, S.: Liouville Quantum Gravity and KPZ, \emph{Inventiones Mathematicae} \textbf{185} (2), 333-393 (2011).


\bibitem{Ry1}
El-Showk, S., Paulos M.F., Poland D., Rychkov S., Simmons-Duffin  D.,  Vichi A.: Solving the 3D Ising model with the conformal bootstrap, \emph{Phys. Rev. D} {\bf 86}, 025022 (2012).

\bibitem{Ry2}
El-Showk, S., Paulos M.F., Poland D., Rychkov S., Simmons-Duffin  D.,  Vichi A.: Solving the 3D Ising model with the conformal bootstrap II. c-minimization and precise critical exponents , \emph{ Phys. Rev. D} {\bf 86}, 025022 (2012).



\bibitem{fybu}
Fyodorov Y., Bouchaud J.-P.: Freezing and extreme value statistics in a Random
Energy Model with logarithmically correlated potential, \emph{J. Phys.A: Math.Theor} {\bf 41}, 372001 (2008).
 

 \bibitem{GouLi} 
 Goulian M, Li M.: Correlation functions in Liouville theory, \emph{Phys. Rev. Lett.} {\bf 66}, 2051 (1991).
 
 
 
 \bibitem{HRV}
 Huang Y., Rhodes R., Vargas V.: Liouville Quantum Gravity on the unit disk, arXiv:1502.04343. 
  
 
 \bibitem{KaraSh}
 Karatzas I. Shreve S.: Brownian motion and stochastic calculus, Springer-Verlag.
 
\bibitem{cf:KPZ} Knizhnik, V.G., Polyakov, A.M., Zamolodchikov, A.B.: Fractal structure of 2D-quantum gravity, \emph{Modern Phys. Lett A}, \textbf{3}(8), 819-826 (1988).
 
 
 
 \bibitem{KoPe}
 Kostov I.K., Petkova V.B.: Bulk correlation functions in 2D quantum gravity, \emph{Theoretical and mathematical physics} {\bf 146} (1), 108-118 (2006).  
 
  \bibitem{cargese} Kupiainen A., Constructive Liouville Conformal Field Theory,  \href{https://arxiv.org/abs/1611.05243}{arXiv:1611.05243}.

 \bibitem{KRV}
 Kupiainen A., Rhodes R., Vargas V.: Local conformal structure of Liouville Quantum Gravity, arXiv:1512.01802. 

\bibitem{KRV1}
 Kupiainen A., Rhodes R., Vargas V.: Integrability of Liouville theory: proof of the DOZZ Formula, arXiv:1707.08785.  
 



\bibitem{nakayama}
Nakayama Y.: Liouville field theory: a decade after the revolution, \emph{Int.J.Mod.Phys. A} \textbf{19}, 2771-2930 (2004).
\bibitem{OPS}
 O'Raifeartaigh L., Pawlowski J.M., Sreedhar V.V.: The Two-exponential Liouville Theory and the Uniqueness of the Three-point Function, \emph{Physics Letters B} {\bf 481} (2-4), 436-444 (2000).
 
 
\bibitem{Pol}
Polyakov A.M.: Quantum geometry of bosonic strings, \emph{Phys. Lett. } \textbf{103B} 207 (1981).

\bibitem{Pol1}
Polyakov A.M.: From Quarks to Strings, \href{https://arxiv.org/abs/0812.0183}{arXiv:0812.0183}. 

\bibitem{remy}
Remy G.: The Fyodorov-Bouchaud formula and Liouville Conformal Field theory, \href{https://arxiv.org/abs/1710.06897}{arXiv:1710.06897}. 
 
\bibitem{Rib}
Ribault S.: Conformal Field theory on the plane, \href{https://arxiv.org/abs/1406.4290}{arXiv:1406.4290}. 

\bibitem{Rib1}
Ribault S., Santachiara R.: Liouville theory with a central charge less than one, \emph{Journal of High Energy physics} {\bf 8}, 109 (2015). 

\bibitem{Rnew10} Rhodes, R. Vargas, V.: KPZ formula for log-infinitely divisible multifractal random measures,  \emph{ESAIM Probability and Statistics} \textbf{15}, 358-371 (2011).

\bibitem{review} 
Rhodes R., Vargas, V.: Gaussian multiplicative chaos and applications: a review,   \emph{Probability Surveys} {\bf 11}, 315-392 (2014).

\bibitem{RV}
Rhodes R., Vargas V.: Lecture notes on Gaussian multiplicative chaos and Liouville Quantum Gravity, \href{https://arxiv.org/abs/1602.07323}{arXiv:1602.07323}.  


\bibitem{Scho}
Schomerus V.: Rolling Tachyons from Liouville theory, \emph{Journal of High Energy Physics} {\bf 11} (2013). 

\bibitem{seiberg}
Seiberg N.: Notes on Quantum Liouville Theory and Quantum  Gravity, \emph{Progress of Theoretical Physics}, suppl. 102 (1990).

\bibitem{She07}
Sheffield S.: Gaussian free fields for mathematicians, \emph{Probab. Th. Rel. Fields} \textbf{139}, 
521-541 (2007).


\bibitem{TT} Takhtajan L., Teo L.-P.: Quantum Liouville Theory in the Background Field Formalism I. Compact Riemann Surfaces, \emph{Communications in mathematical physics} {\bf 268} (1), 135-197 (2006).  





\bibitem{Tesc} Teschner J.: On the Liouville three point function, \emph{Phys. Lett. B} {\bf 363}, 65-70 (1995).

\bibitem{Tesc1} Teschner J.: Liouville Theory Revisited, \emph{Class.Quant.Grav.} {\bf 18}, R153-R222 (2001).

\bibitem{Var} 
Vargas V.: IHES lectures available at https://www.youtube.com/watch?v=BU8VO6ps59s.


\bibitem{Williams}
Williams, D.: Path Decomposition and Continuity of Local Time for One-Dimensional Diffusions, I, \emph{Proceedings of the London Mathematical Society} {\bf s3-28} (4), 738-768 (1974).


\bibitem{Za}
Zamolodchikov A.B.: Three-point function in the minimal Liouville gravity, \emph{Theoretical and Mathematical Physics} {\bf 142} (2), 183-196 (2005).  



\bibitem{ZaZaarxiv}
Zamolodchikov A.B., Zamolodchikov A.B.:  Structure constants and conformal bootstrap
in Liouville field theory, arXiv:hep-th/9506136.

\bibitem{ZaZa}
Zamolodchikov A.B., Zamolodchikov A.B.: Conformal bootstrap
in Liouville field theory, \emph{ Nuclear Physics B}
{\bf 477} (2), 577-605 (1996).



\end{thebibliography}
\end{document}